\documentclass[11pt,twoside]{article}
\usepackage{amsfonts}
\usepackage{fancyhdr}
\usepackage{titlesec}
\usepackage{cite}
\usepackage{ifthen}
\usepackage{amssymb}
\usepackage{fancyhdr}
\usepackage{titlesec}
\usepackage[arrow,matrix]{xy}
\usepackage{pifont}
\usepackage{stmaryrd}
\usepackage{setspace}
\usepackage{indentfirst}
\usepackage{amsmath,amssymb,amscd,bbm,amsthm,mathrsfs,dsfont}
\usepackage{graphicx}
\input amssym.def

\def\setzero{\setcounter{equation}{0}}

\titleformat{\section}{\centering\large\bfseries}{\S\arabic{section}}{1em}{}
\newboolean{first}
\setboolean{first}{true}

\begin{document}

\setlength\abovedisplayskip{2pt}
\setlength\abovedisplayshortskip{0pt}
\setlength\belowdisplayskip{2pt}
\setlength\belowdisplayshortskip{0pt}

\begin{center}
{\Large \bf Optimal stochastic control and optimal consumption and portfolio  with G-Brownian motion
\footnote{Corresponding author: Weiyin Fei, email: wyfei@dhu.edu.cn.}}

\vskip12pt

{\rm Weiyin Fei \quad Chen Fei}

\vskip12pt
\small{ School of Mathematics and Physics, Anhui Polytechnic University, Wuhu, Anhui, P.R. China, 241000\\}

\end{center}

\vskip12pt

\noindent{\bf Abstract:} By the calculus of Peng's G-sublinear expectation and G-Brownian motion on a sublinear expectation space $(\Omega, {\cal H}, \hat{\mathbb{E}})$, we first set up an optimality principle of stochastic control problem. Then we investigate an optimal consumption and portfolio decision with a volatility ambiguity by the derived verification theorem. Next the two-fund separation theorem is explicitly obtained. And an illustrative example is provided.

  \vskip12pt


 \vskip12pt

 \noindent{\bf MR(2000) Subject Classification} 60H30; 93E20; 91G10

\section{Introduction} \setzero

In the real world, we are often faced with two kinds of uncertainties, i.e., probabilistic uncertainty and Knightian uncertainty (model uncertainty or ambiguity) (c.f. Ellsberg \cite{El} and Knight \cite{Kn}). Knightian uncertainty is due to incomplete information, vague data, imprecise probability etc. Hence, uncertainty is inherent in most real-world systems, it placed many disadvantages or sometimes, surprisingly, advantages on humankind's efforts, which are usually associated with the request for optimal results.

First, let us recall simply the history of development of the Brownian motion theory since Brownian motion is often believed the source of noise of an actual system.
In fact, in 1827, the botanist Robert Brown, looking through a microscope at particles found in pollen grains in water, noted that the particles moved through the water but was not able to determine the mechanisms that caused this motion. Atoms and molecules had long been theorized as the constituents of matter, and many decades later, Albert Einstein \cite{E} published a paper in 1905 that explained in precise detail how the motion that Brown had observed was a result of the pollen being moved by individual water molecules. The mathematical model of Brownian motion has numerous real-world applications. For instance, stock market fluctuations was early cited by Bachelier \cite{Ba} who believed that the stock's prices follows Brownian motion. His thesis, which discussed the use of Brownian motion to evaluate stock options, is historically the first paper to use advanced mathematics in the study of finance. Hence, Bachelier is considered a pioneer in the study of financial mathematics and stochastic processes.
Brownian motion is among the simplest of the continuous-time stochastic (or probabilistic) processes, and it is a limit of both simpler and more complicated stochastic processes (see random walk and Donsker's theorem). This universality is closely related to the universality of the normal distribution. In both cases, it is often mathematical convenience rather than the accuracy of the models that motivates their use. This is because Brownian motion, whose time derivative is everywhere infinite, is an idealised approximation to actual random physical processes, which always have a finite time. However, by the way, in many actual systems, it might be believed that the ``source of noise'' might not be a Brownian motion, just be a fractional Brownian motion which admits self-similarity and lang-range dependence (c.f. Mandelbrot and Van Ness \cite{MN}). Moreover, the related work for a fractional Brownian motion is proceeded, such as Fei et al. \cite{FXZ}, and references therein.

Second, it is necessary to study how we make a decision when we are often faced with probabilistic uncertainty and Knightian uncertainty. To this end, many researchers investigate the characteristics of model uncertainty in order to provide a framework for theory and applications. Choquet \cite{Ch} put forward to the notion of capacity which is a nonadditive measure in 1953. In 2006, Peng \cite{P2006} think that a classical Brownian motion has not characterized ambiguous volatility, hence he put forward G-Brownian motion and the related It\^o calculus which started a new area of research.  Recently, the theory of G-expectation  and G-Brownian motion on sublinear expectation space, which  provides the new perspective for the stochastic calculus under Knightian uncertainty, is further discussed in Peng and his coauthors \cite{P2008,P2010,HP,LP}. In Soner et al. \cite{STZ}, the martingale representation
theorem under G-expectation is proven.
The recent developments on problems of probability
model with ambiguity by using the Peng's notion of nonlinear expectations and, in
particular, sublinear expectations show that a nonlinear expectation is a monotone and constant preserving functional defined on a linear space of random variables.
A sublinear expectation can be represented as the upper expectation of
a subset of linear expectations. In
most cases, this subset is often treated as an uncertain model of probabilities ( for a financial decision maker, this multiple priors set reflects the uncertainty degree of the decision maker).  The sublinear expectation theory provides many rich, flexible and elegant tools.
We emphasize the term ``expectation''
rather than the well-accepted classical notion ``probability'' and its non-additive
counterpart ``capacity''. The notion of expectation has its direct meaning of ``mean'' which is not necessary to be derived from the corresponding ``relative frequency'' which is the origin of the probability measure.

Third, since the study of optimal stochastic controls from the theoretical point of view as well as for their applications (including financial economics) is important, our question is how we provide a framework of a stochastic control system with ambiguity. In fact, we know that the modern optimal control theory began its life at the end of World War II. The main starting point seems to be (differential) games. Such a study created a perfect environment for the birth of Bellman's dynamic programming method (see Bellman \cite{Be}). According to Pontryagin \cite{Pon}, with him and his coauthors' extensive study, Pontryagin's maximum principle was announced in 1956.  In a classical stochastic control theory, the method of dynamic programming is a mathematical technique for making a sequence of interrelated decisions, which can be applied to many optimization problems (including optimal stochastic control). The basic idea of this method applied to optimal stochastic controls is to consider a family of optimal stochastic control problems with different initial times and states, to establish relations among these problems by the Hamilton-Jacobi-Bellman (HJB) equation, which is second-order partial differential equation. If the HJB equation is solvable, then one can obtain an optimal stochastic control via taking the maximizer/minimizer of the Hamiltonian in HJB equation. The applications of duality to optimal stochastic control are given in Bismit \cite{B}. The controlled Markov processes and viscosity solutions were explored by Fleming and Soner \cite{FS}. The elegant investigation of the classical stochastic control is provided by Yong and Zhou \cite{YZ}.
 Since the classical stochastic control cannot consider a model uncertainty (especially, volatility ambiguity), it is necessary to investigate a system in some complex environment with ambiguity by a calculus of sublinear expectation. We know that there exists a large literature analyzes stochastic differential systems with Wiener's noise (see, e.g., Karaztas and Shreve \cite{KS1991}, Mao \cite{Mao}, $\O$ksendal \cite{O}, Revus and Yor \cite{RY}). Recently, under uncertainty, a kind of stochastic differential equation driven by G-Brownian motion is studied by Gao \cite{G} and Peng \cite{P2010}. To our best knowledge, the study of the optimal stochastic control based on Peng's sublinear expectation theory is not still found. So in this paper we will investigate the dynamic programming principle and  the corresponding  HJB equations.

Next, in the investigation of a financial market, we formulate  a  model  of  continuous  time  utility  with ambiguity or model uncertainty, where                     the  multiple  priors  (or  maxmin)  model  of  preference
put forth by Gilboa and Schmeidler \cite{GS} for a static setting is adapted.  Under the model uncertainty, by using the method of the robust control Hansen et al. \cite{HS} discussed min-max expected utility where an ambiguous volatility is not considered. Also, in Chen and Epstein \cite{CE}, Fei \cite{FeiINS,FeiSM,FeiL}, all probability measures entertained by a decision-maker are assumed to be equivalent to a
fixed reference probability measure. So there is no ambiguity about which scenarios are possible.
We know that from an
economics perspective, the assumption of equivalence seems far from innocuous.
Informally, if her environment is complex, how could the decision-maker come to
be certain of which scenarios regarding future asset prices and rates of return are possible?  Especially, ambiguity about volatility implies ambiguity about which scenarios are possible, at least in a continuous time setting.
    Since volatility is not directly observable and
must be inferred from observation such as realized asset returns and prices, we are led to develop a model of preference that accommodates ambiguity about
volatility. Epstein and Ji \cite{EJ1,EJ2} generalized the Chen-Epstein model and maintained a separation between risk aversion
and intertemporal substitution, where the asset pricing with volatility ambiguity is explored. In order to apply the model with Knightian uncertainty in the  continuous  time  setting, we will study the optimal consumption and portfolio model of an agent with volatility ambiguity by the obtained results of an optimal stochastic control. On a study of optimal consumption and portfolio in the continuous time setting, we can date back to the Merton's work in 1969 and 1971 (c.f. Merton \cite{M1969,M1971}). Later, a large literature investigates this subject with all kinds of conditions, for instance, Karaztas and Shreve \cite{KS1998}, Duffie \cite{Du}, Fei et al. \cite{FeiST, FW2002,FW2003}, and references therein.

In a word, when we are faced with Knightian uncertainty, the stochastic control systems perturbed by G-Brownian motion will be  important for characterizing the real world with both randomness and ambiguity. Specifically, it is necessary to study the problem of optimal stochastic controls with ambiguity in a similar manner as in classical ones. To this end, the organization of the paper is as follows. In Section 2, we give preliminaries and provide some lemmas.
In Section 3, the principles of optimality of optimal stochastic controls with G-Brownian noise are proven. Section 4 studies the optimal consumption and portfolio as a financial application of optimal stochastic control with ambiguity, where the optimal consumption and portfolio decision is obtained through HJB equation, and an illustrative example is given.  Finally, Section 5 concludes.

\section{Preliminaries}\setzero

 In this section, we first give the notion of sublinear expectation space $(\Omega, {\cal H}, {\mathbb{E}})$, where $\Omega$ is a given state set and $\cal H$ a linear space of real valued functions defined on $\Omega$. The space $\cal H$ can be considered as the space of random variables. The following concepts come from  Peng \cite{P2010}.

\vskip12pt

{\noindent\bf Definition 2.1}  A sublinear expectation $\mathbb E$ is a functional $\mathbb  E$: ${\cal H}\rightarrow {\mathbb  R}$ satisfying

\noindent(i) Monotonicity: ${\mathbb{E}}[X]\geq {\mathbb{E}}[Y]$ if $X\geq Y$;

\noindent(ii) Constant preserving: ${\mathbb{E}}[c]=c$;

\noindent (iii) Sub-additivity: For each $X,Y\in {\cal H}$, ${\mathbb{E}}[X+Y]\leq {\mathbb{E}}[X]+{\mathbb{E}}[Y]$;

\noindent(iv) Positivity homogeneity: ${\mathbb{E}}[\lambda X]=\lambda{\mathbb{E}}[X]$ for $\lambda\geq 0$.

\vskip12pt
We denote by $\mathbb{S}(d)$ the collection of all $d\times d$ symmetric matrices (resp., ${\mathbb{S}^{>0}(d)}$ the space of all $d\times d$ positive-definite matrices). $\Sigma\subset{\mathbb{S}^{>0}(d)}$ is the bounded, convex and closed subset which is of form $\Sigma=\{\Lambda\in {\mathbb S}^{>0}: \underline{\sigma}^2I_d\leq \Lambda\leq \bar{\sigma}^2 I_d\} $ for some positive constants $\underline{\sigma}^2$ and $\bar{\sigma}^2$.

\vskip12pt

{\noindent\bf Definition 2.2} Let $(\Omega,{\cal H}, {\mathbb{E}})$ be a sublinear expectation space. $(X_t)_{t\geq0}$ is called a $d$-dimensional stochastic process if for each $t\geq0$, $X_t$ is a $d$-dimensional random vector in $\cal H$.

A $d$-dimensional process $(B(t))_{t\geq0}$ on a sublinear expectation space $(\Omega,{\cal H}, {\mathbb{E}})$ is called a G-Brownian motion if the following properties are satisfies:

\noindent(i) $B_0(\omega)=0$;

\noindent (ii) for each $t,s\geq 0$, the increment $B(t+s)-B(t)$ is $N(\{0\}\times s\Sigma)$-distributed and independent from $(B_{t_1},B_{t_2},\cdots, B_{t_n})$, for each $n\in {\mathbb N}$ and $0\leq t_1\leq \cdots\leq t_n\leq t$.

\vskip12pt

Next  we denote by $\Omega =C_0^d({\mathbb R}_+)$ the space of all ${\mathbb R}^d$-valued continuous paths $(\omega_s)_{s\geq0}$, with $\omega_0=0$, equipped with the distance
$$\rho(\omega^1,\omega^2):=\sum\limits_{i=1}^\infty2^{-i}[(\max\limits_{t\in[0,i]}|\omega^1_t-\omega^2_t|)\wedge 1].$$
Define $\Omega_s^t:=\{\omega_{\cdot\wedge s}-\omega_{\cdot\wedge t}: \omega\in \Omega\}$.
More details on the notions of G-expectation ${\hat{\mathbb E}}:L_{ip}(\Omega)\rightarrow{\mathbb R}$ and G-Brownian motion on the sublinear expectation space $(\Omega,L_{ip}(\Omega), {\hat{\mathbb{E}}})$ are referred to Peng \cite{P2010}. Moreover, from Denis et al. \cite{DHP} and Theorem VI-2.5 of Peng \cite{P2010}, we know that there exists a weak compact family of probability measures ${\cal P}$ on $(\Omega, {\cal B}(\Omega))$ such that $\hat{\mathbb E}[X]=\max_{P\in{\cal P}}E_P[X],~\forall X\in ~L_{ip}(\Omega)$, where $E_P[\cdot]$ is a linear expectation with respect to $P$. Thus a property holds ``quasi-surely" (q.s.) if it holds almost surely for each prior $P\in {\cal P}$.

We now give the definition of It\^o integral. For simplicity, here we only introduce It\^o integral with respect to 1-dimensional G-Brownian motion with $G(\alpha):=\frac{1}{2}\hat{\mathbb{E}}[\alpha B_1^2]=\frac{1}{2}(\bar{\sigma}^2\alpha^+-\underline{\sigma}^2\alpha^-)$, where $\hat{\mathbb{E}}[B_1^2]=\bar{\sigma}^2, -\hat{\mathbb{E}}[- B_1^2]=\underline{\sigma}^2$,
$0<\underline{\sigma}\leq \bar{\sigma}<\infty$.

 Let $L^p_G(\Omega), p\geq1$ be the completion of $L_{ip}(\Omega)$ under the norm $\|X\|_p:=(\hat{\mathbb E}|X|^p)^{1/p}$. We consider the following type of simple processes: for a given partition $\pi_T=(t_0,\cdots,t_N)$
of $[0, T]$ we get that
$$\eta_t(\omega)=\sum_{k=0}^{N-1}\xi_k(\omega)I_{[t_k,t_{k+1})}(t),$$
where $\xi_k\in L^p_G(\Omega_{t_k}), k=0,1,\cdots,N-1$ are given. The collection of these processes is denoted by $M_G^{p,0}(0,T)$. We denote by $M_G^p(0,T)$ the completion of $M_G^{p,0}(0,T)$  with the norm
 $$\|\eta\|_{M_G^p(0,T)}:=\left\{\hat{\mathbb{E}}[\int_0^T|\eta_t|^pdt]\right\}^{1/p}<\infty.$$
For $\eta\in M_G^p(0,T)$, the definition of stochastic integral $\int_0^T\eta_t dB(t)$ is referred to Peng \cite{P2010}.

We now consider the following G-stochastic differential equation (G-SDE)
$$X_t^\nu=X_0^\nu+\int_0^t\alpha_s^\nu ds+\sum_{i,j=1}^d\int_0^t\eta_s^{\nu i j}d<B^i,B^j>_s+\sum\limits_{j=1}^d\int_0^t\beta_s^{\nu j}dB^j(s),~\nu=1,\cdots n.$$
For convenience, we generalize G-It\^o formula from Peng \cite{P2010} and Gao \cite{G} as follows.

 \vskip12pt

  {\noindent \bf Lemma 2.3} (G-It\^o Formula) {\sl Let $\Psi$ be a $C^{1,2}$-function on $\mathbb{R}_+\times\mathbb{R}^n$ such that $\Psi_{x^\mu x^\nu}(t,x)$ is local Lipschitz in $x$, i.e., for each $m\geq1$ there exists a constant $L_m>0$ such that for any $x,y\in {\mathbb R}^n$ with $|x|\leq m, |y|\leq m$,
  $$\max\limits_{1\leq\mu,\nu\leq n}|\Psi_{x^\mu x^\nu}(t,x)-\Psi_{x^\mu x^\nu}(t,y)|<L_m,~\forall t\in {\mathbb R}_+.$$
  If for each process $\xi$ in
  $$\cup_{\mu,\nu=1,\cdots,n;i,j=1,\cdots,d}\{\alpha^{\nu j},\beta^{\nu j},\eta^{\nu i j},\Psi_t,\Psi_{x^\nu}\alpha^{\nu j},\Psi_{x^\nu}\beta^{\nu j},\Psi_{x^\nu}\eta^{\nu ij},\Psi_{x^\mu x^{\nu}}\beta^{\mu i}\beta^{\nu j}\},$$
  for each $ t\in[0,T],$
  $$\sup\limits_{t\in[0,T]}\hat{\mathbb E}[|\xi_t|^2]<\infty,~\xi_t\in L_G^2(\Omega_t),~\lim\limits_{\delta\rightarrow0}\sup\limits_{s:|s-t|\leq \delta}\hat{\mathbb E}[|\xi_s-\xi_t|^2]=0, $$
  then $\{\Psi_t,\Psi_{x^\nu}\alpha^{\nu j},\Psi_{x^{\nu}}\beta^{\nu j}, \Psi_{x^\nu}\eta^{\nu i j},\Psi_{x^{\nu}x^{\mu}}\beta^{\nu i}\beta^{\nu j}\}\subset M_G^2(0,T)$, and for each $t\in [0,T]$ we have, in $L_G^2(\Omega_t)$,
   $$
   \begin{array}{ll}
   &\Psi(t,X_t)-\Psi(s,X_s)\\
   &=\sum\limits_{\nu=1}^n\sum\limits_{j=1}^d\int_s^t\Psi_{x^\nu}(r,X_r)\beta_r^{\nu j}dB^j(r)+\sum\limits_{\nu=1}^n\int_s^t[\Psi_t (r,X_r)+\Psi_{x^\nu}(r,X_r)\alpha_r^\nu ]dr
   \end{array}$$
    $$
   +\sum\limits_{\mu,\nu=1}^n\sum\limits_{i,j=1}^d\int_s^t[\Psi_{x^\nu}(r,X_r)\eta_r^{\nu i j}+\frac{1}{2}\Psi_{x^\nu x^\mu}(r, X_r)\beta_r^{\mu i}\beta^{\nu j}_r]d<B^i,B^j>_r.
$$
   
   }

  \vskip12pt
  Next, from Proposition III-2.8 in Peng \cite{P2010}, we have the following result.

  \vskip12pt.

  {\bf \noindent Lemma 2.4} {\sl Let $X,Y\in L_G^1(\Omega)$ such that $\hat{\mathbb E}[Y|\Omega_t]=-\hat{\mathbb E}[-Y|\Omega_t]$, for some $t\in[0,T]$. Then we have
  $$
  \begin{array}{rl}
  \hat{\mathbb E}[X+Y|\Omega_t]=&\hat{\mathbb E}[X|\Omega_t]+\hat{\mathbb E}[Y|\Omega_t],\\
 -\hat{\mathbb E}[-(X-Y)|\Omega_t]=&-\hat{\mathbb E}[-X|\Omega_t]-(-\hat{\mathbb E}[-Y|\Omega_t]).
  \end{array}
  $$
  }

  \vskip12pt

    Let $\prec p, q\succ$ denote the inner product of  vectors $p,q\in \mathbb{R}^d$. We define
   $$\prec A, B\succ:=tr({A^\top B})=\sum_{i=1}^l\sum_{j=1}^da_{ij}b_{ij},~ \mbox{for}~A,B\in\mathbb{R}^{l\times d},$$
   and
   $$G(A):=\frac{1}{2}\sup\limits_{\Lambda\in \Sigma}\prec A, \Lambda\succ,~A\in \mathbb{S}(d),$$
   $$\widetilde{G}(A):=\frac{1}{2}\inf\limits_{\Lambda\in \Sigma}\prec A, \Lambda\succ,~A\in \mathbb{S}(d).$$
From Corollary III-5.7 in Peng \cite{P2010}, we know $<B>_t\in t\Sigma:=\{t\times\gamma:\gamma\in\Sigma\}$.

From F\"ollmer \cite{Fo}, Karandikar \cite{Ka}, Soner et al. \cite{STZ}, and Epstein and Ji \cite{EJ2}, there exists a probability measure $P^{(\varrho_t)}\in {\cal P}$ with $v_t:=\varrho_t\varrho^\top_t\in \Sigma$  such that
$$d<B>_t=v_t dt=\varrho_t\varrho_t^\top dt, ~dt\times P^{(\varrho_t)}\mbox{-a.e.},$$
which shows that
$$\hat{\mathbb E}[\prec A,v_t\succ]=2G(A),~-\hat{\mathbb E}[-\prec A,v_t\succ]=2\widetilde{G}(A),~ \forall A\in {\mathbb S}(d). \eqno(2.1)$$

\vskip12pt

\section{Stochastic principle of optimality}\setzero

Let $\bf U$ be a separable metric space. For a control process $u(\cdot):[0,T]\times \Omega\rightarrow {\bf U}$, our controlled system follows
$$
\begin{array}{ll}
\left\{\begin{array}{ll} dx(t)=f(t, x(t), u(t))dt+g(t,
x(t),u(t))dB(t), \quad
t\in[0,T],\\
x(0)=x,
\end{array}
\right.
\end{array}
\eqno(3.1)
$$
where $B(\cdot)$ is a G-Brownian motion, and
$$
\begin{array}{ll}
&f=(f^1,\cdots,f^m)^\top:[0,T]\times{\mathbb{R}}^m\times
{\bf U}\rightarrow {\mathbb R}^m,\\
& g=(g^{kl})_{m\times d}: [0,T]\times{\mathbb{R}}^m\times{\bf U}\rightarrow {\mathbb R}^{m\times d}.\\
\end{array}
$$

 First of all, we consider that a controller makes a decision with conservative attitude on Knightian uncertainty.  Thus, we may suppose that the cost functional of our control problem with multiple-priors ${\cal P}$ is as
follows
$$
\begin{array}{rl}
J(x; u(\cdot)):=&\max\limits_{P\in{\cal P}}E_P\left[\int_0^T\Phi_1(t, x(t), u(t))dt+\Phi_2(x(T))\right]\\
=&\mathbb{\hat{E}}\left[\int_0^T\Phi_1(t, x(t), u(t))dt+\Phi_2(x(T))\right].
\end{array}
 \eqno(3.2)
$$
We now introduce an assumption.

\vskip12pt

\noindent{\bf Assumption 3.1} The maps $f:[0,T]\times {\mathbb R}^m\times {\bf U}\rightarrow {\mathbb R}^m, g:[0,T]\times{\mathbb R}^m\times{\bf U}\rightarrow{\mathbb R}^{m\times d}, \Phi_1: [0,T]\times{\mathbb R}^m\times{\bf U}\rightarrow{\mathbb R}$, and $\Phi_2: {\mathbb R}^m\rightarrow{\mathbb R}$
are uniformly continuous, and there exists a constant $L>0$ such that for $\phi(t,x,u)=f(t,x,u),g(t,x,u),\Phi_1(t,x,u),\Phi_2(x)$,
$$
\begin{array}{ll}
\left\{\begin{array}{ll} |\phi(t,x,u)-\phi(t,\bar{x},u)|\leq L|x-\bar{x}|,~\forall t\in[0,T],x,\bar{x}\in{\mathbb R}^m,u\in{\bf U},\\
\phi(t,0,u)|\leq L,~\forall(t,u)\in[0,T]\times {\bf U}.
\end{array}
\right.
\end{array}
$$

\vskip12pt

 Subsequently, we need to consider a family of optimal control problems
with different initial times other than zero and states along a given state
trajectory in order to apply the dynamic programming technique. To this end, we give the following dynamics. For any fixed $(t,
x)\in[0,T]\times {\mathbb{R}}^m$, consider the state equation
$$
\begin{array}{ll}
\left\{\begin{array}{ll} dx(s)=f(s, x(s), u(s))ds+g(s,
x(s),u(s))dB(s),~\forall s\in[t,T],\\
x(t)=x.
\end{array}
\right.
\end{array}
\eqno(3.3)
$$

\vskip12pt

{\noindent \bf Problem 3.2} A controller's cost functional with conservative attitude on Knightian uncertainty is
$$
J(t,x; u(\cdot))
=\mathbb{\hat{E}}\left[\int_t^T\Phi_1(s, x(s), u(s))ds+\Phi_2(x(T))\right],
 $$
 where $u(\cdot)$ is a control process. Let us denote the value function by
 $$V(t, x):=\inf_{u(\cdot)\in {\cal U}[t,T]}J(t,
x;u(\cdot)),$$
 where ${\cal U}[t,T]$ denotes the set of all $u$ satisfying

\noindent(i) $u:[t,T]\times \Omega\rightarrow {\bf U}$ is an $\{\Omega_s^t\}$-adapted process on a sublinear expectation space $(\Omega,{\cal H},\hat{\mathbb E})$.

\noindent(ii) Under a control $u(\cdot)$, for each $x\in{\mathbb R}^m$ equation (3.3) admits a unique solution $x(\cdot)$  on a sublinear expectation space $(\Omega,{\cal H},\hat{\mathbb E})$.

\noindent(iii) $\Phi_1(\cdot,x(\cdot),u(\cdot))\in M^1_G(0,T)$  and $\Phi_2(x(T))\in L_G^1(\Omega_T)$.

\vskip12pt

{\noindent \bf Rremark 3.3}  It is known (cf. Peng \cite{P2010}) that under
Assumption 3.1, Eq. (3.1) admits the unique
continuous solutions $x(t)=x(t,u(t))$ on $t\in[0,T]$ for each $u(\cdot)\in {\cal U}[0,T]$.

\vskip12pt

For our aim, we now provide the following lemma.

\vskip12pt

{\noindent \bf Lemma 3.4} {\sl Let Assumption 3.1 hold.  Then,  the equation (3.3) with $x(0)=x$ admits a unique solution $x(\cdot)$ such that for any $T>0$ and $\ell\geq 1$,
$$\hat{\mathbb E}\max\limits_{0\leq s\leq T}|x(s)|^\ell\leq K_{\ell,T}(1+|x|^\ell)\eqno(3.4)$$
and
$$\mathbb{\hat{E}}|x(t)-x(s)|^\ell\leq K_{\ell,T}(1+|x|^\ell)|t-s|^{\ell/2}, ~\forall s,t\in[0,T].\eqno(3.5)$$
Moreover, if $\hat{x}(0)=y\in {\mathbb R}^m$ and $\hat{x}(\cdot)$ is the corresponding solution of (3.3), then for any $T>0$, there exists a $K_{\ell,T}>0$ such that

$$\mathbb{\hat{E}}\max\limits_{0\leq s\leq T}|x(s)-\hat{x}(s)|^\ell\leq K_{\ell,T}|x-y|^\ell.\eqno(3.6)$$
}

\vskip12pt

\begin{proof} Since a controller has a multiple-priors set $\cal P$ which is weakly compact, for each $P\in {\cal P}$, there exists a $K_{\ell,T}(P)>0$ depending $P$ such that (3.4)-(3.6) hold in a way similar to the proof of Theorem I-6.16 in \cite{YZ}.
Set $K_{\ell,T}=\max_{P\in{\cal P}} K_{\ell,T}(P)$. Thanks to $\cal P$ being weakly compact, we have that $K_{\ell,T}<\infty$. Indeed, if the above claim is false, then for each positive integer $N>0$ there exists a $P_N\in {\cal P}$ such that $K_{\ell,T}(P_N)>N$ satisfying $K_{\ell,T}(P_{N})\uparrow\infty$. On the other hand, by the weak compactness of ${\cal P}$ there exists a subsequence $\{P_{N_k}\}$ of $\{P_N\}$ with $K_{\ell,T}(P_{N_k})\uparrow\infty$ such that $P_{N_k}$ converge weakly to $P\in {\cal P}$ with $K_{\ell,T}(P)<\infty$ which shows $K_{\ell,T}(P_{N_k})\rightarrow K_{\ell,T}(P)$. This is  a contraction. Thus the proof is complete.\end{proof}

\vskip12pt

{\noindent \bf Proposition 3.5} {\sl Let Assumption 3.1 hold. Then for some constant $K>0$, the value function $V(t,x)$ satisfies
$$|V(t,x)|\leq K(1+|x|),~ \forall (t,x)\in [0,T]\times {\mathbb R}^m,~\eqno(3.7)$$
$$
\begin{array}{ll}
|J(t,x;u(\cdot))-J(\bar{t},\bar{x};u(\cdot))|\vee|V(t,x)-V(\bar t, \bar x)|\leq &K(|x-\bar{x}|+(1+|x|\vee|\bar{x}|)|t-\bar{t}|^{1/2}),\\
&\forall t,\bar{t}\in[0,T],x,\bar{x}\in{\mathbb R}^m, u(\cdot)\in{\cal U}[t,T].
\end{array}
\eqno(3.8)
$$
}

\vskip12pt

\begin{proof} Fix $(t,x)\in[0,T]\times {\mathbb R}^m$. For each $u(\cdot)\in {\cal U}[t,T]$, by Lemma 3.4 we have
$${\mathbb E}\sup\limits_{s\in[t,T]}|x(s)|\leq K(1+|x|). \eqno(3.9)$$
Thus, in terms of Assumption 3.1 and (3.9), we have
$$|J(t,x;u(\cdot))|\leq K(1+|x|),~\forall u(\cdot)\in {\cal U}[t,T],$$
which shows (3.7) through taking the infimum in $u(\cdot)\in{\cal U}[t,T]$.

Next we let $0\leq t\leq\bar{t}\leq T$ and $x,\bar{x}\in {\mathbb R}^m$. For any control $u(\cdot)\in {\cal U}[t,T]$, let $x(\cdot)$ and $\bar{x}(\cdot)$ be the states corresponding to $(t,x,u(\cdot))$ and $(\bar{t},\bar{x},u(\cdot))$, respectively. Therefore, by Lemma 3.4, we obtain
$$
\mathbb{\hat{E}}\sup\limits_{s\in[\bar{t},T]}|x(s)-\bar{x}(s)|\leq K (|x-\bar{x}|+(1+|x|\vee|\bar{x}|)|t-\bar{t}|^{1/2}).
$$
Thus, by Assumption 3.1 we get
$$
|J(t,x;u(\cdot))-J(\bar{t},\bar{x};u(\cdot))|\leq  K (|x-\bar{x}|+(1+|x|\vee|\bar{x}|)|t-\bar{t}|^{1/2}).\eqno(3.10)
$$
Taking the infimum in $u(\cdot)\in{\cal U}[t,T]$, together with (3.10), we obtain (3.8). Thus the proof is complete.
\end{proof}

\vskip12pt

Next we give the dynamic programming principle.

\vskip12pt

{\noindent \bf Theorem 3.6} (Principle of optimality I) {\sl Let Assumption 3.1 hold. Then for any $(t,x)\in [0,T]\times {\mathbb R}^m,$ we have
$$
\begin{array}{ll}
V(t,x)=\inf\limits_{u(\cdot)\in{\cal U}[t,T]}&\mathbb{\hat{E}}\{\int_t^{\bar{t}}\Phi_1(s,x(s;t,x,u(\cdot)),u(s))ds\\
&\quad+ V(\bar{t},x(\bar{t}; t,x,u(\cdot))) \},~\forall 0\leq t\leq\bar{t}\leq T.
\end{array}\eqno(3.11)
$$
}

\vskip12pt

\begin{proof} For any $\varepsilon>0$, there exists a $u(\cdot)\in {\cal U}[t,T]$ such that
$$
\begin{array}{ll}
&V(t,x)+\varepsilon>J(t,x;u(\cdot))\\
&=\mathbb{\hat{E}}\left\{\int_t^T\Phi_1(s,x(s;t,x,u(
\cdot)),u(s))ds+\Phi_2(x(T;t,x,u(\cdot)))\right\}\\
&=\mathbb{\hat{E}}\left\{\int_t^{\bar{t}}\Phi_1(s,x(s;t,x,u(\cdot)),u(s))ds\right.\\
&\quad\left.+\mathbb{\hat{E}}\left[\int_{\bar{t}}^T\Phi_1(s,x(s;t,x,u(\cdot)),u(s))ds+\Phi_2(x(T;t,x,u(\cdot)))|{\Omega^t_{\bar{t}}}\right]\right\}\\
&=\mathbb{\hat{E}}\left\{\int_t^{\bar{t}}\Phi_1(s,x(s;t,x,u(\cdot)),u(s))ds\right.\\
&\quad\left.+\mathbb{\hat{E}}\left[\int_{\bar{t}}^T\Phi_1(s,x(s;\bar{t},x(\bar{t}),u(\cdot)),u(s))ds+\Phi_2(x(T;\bar{t},x(\bar{t}),u(\cdot))|{\Omega^t_{\bar{t}}}\right]\right\}\\
&=\mathbb{\hat{E}}\left\{\int_t^{\bar{t}}\Phi_1(s,x(s;t,x,u(\cdot)),u(s))ds+J(\bar{t},x(\bar{t};t,x,u(\cdot));u(\cdot))\right\}\\
&\geq\mathbb{\hat{E}}\left\{\int_t^{\bar{t}}\Phi_1(s,x(s;t,x,u(\cdot)),u(s))ds+V(\bar{t},x(\bar{t};t,x,u(\cdot)))\right\},
\end{array}
$$
where we use
$$
\begin{array}{ll}
J(\bar{t},x(\bar{t});u(\cdot))=&\mathbb{\hat{E}}\left\{\int_{\bar{t}}^T\Phi_1(s,x(s;\bar{t},x(\bar{t}),u(\cdot)),u(s))ds\right.\\
&\quad\left.+\Phi_2(x(T;\bar{t},x(\bar{t}),u(\cdot)))|{\Omega^t_{\bar{t}}}\right\}(\omega)~ \mbox{q.s.}
\end{array}
$$
Hence, letting $\varepsilon\downarrow0$ in above inequality we get that
 $$
V(t,x)\geq \inf\limits_{u(\cdot)\in{\cal U}[t,T]}\mathbb{\hat{E}}\left\{\int_t^{\bar{t}}\Phi_1(s,x(s;t,x,u(\cdot)),u(s))ds+V(\bar{t},x(\bar{t};t,x,u(\cdot)))\right\}.
$$
On the other hand, for any $\varepsilon>0$, by Proposition 3.5 and its proof, there exists a $\delta$ such that whenever $|y-\bar{y}|<\delta$,
 $$|J(\bar{t},y,u(\cdot))-J(\bar{t},\bar{y},u(\cdot))|+|V(\bar{t},y)-V(\bar{t},\bar{y})|\leq\varepsilon,~\forall u(\cdot)\in{\cal U}[\bar{t},T].$$
 Let $\{D_j\}_{j\geq 1}$ be a Borel partition of $\mathbb{R}^m$ satisfying $\cup_{j\geq 1}D_j={\mathbb R}^m$ and $D_i\cap D_j=\emptyset$ if $i\neq j$ with diameter diam $(D_j)<\delta$. Select $x_j\in D_j$. For each $j$ there exists $u_j(\cdot)\in {\cal U}[\bar{t},T]$ such that
 $$J(\bar{t},x_j;u_j(\cdot))\leq V(\bar{t},x_j)+\varepsilon.$$
 Thus for any $x\in D_j$, we have
 $$J(\bar{t},x;u_j(\cdot))\leq J(\bar{t},x_j;u_j(\cdot))+\varepsilon\leq V(\bar{t},x_j)+2\varepsilon\leq V(\bar{t},x)+3\varepsilon.\eqno(3.12)$$
 Therefore there exists a continuous function $\psi_j:[\bar{t}, T]\times C([0,T];{\mathbb R}^d)\rightarrow {\bf U}$ such that
 $$u_j(t,\omega)=\psi_j(t, B(\cdot\wedge t,\omega)-B(\cdot\wedge \bar{t},\omega))~\mbox{q.s.},~\forall t\in [\bar{t},T].$$
 Now, given $u(\cdot)\in {\cal U}[t,T]$, let $x(\cdot)\equiv x(\cdot;t,x,u(\cdot))$ denote the corresponding state trajectory for Problem 3.2. Define a new control
 $$
\begin{array}{ll}
\widetilde{u}(t,\omega)=\left\{\begin{array}{ll}u(t,\omega),&\mbox {if}~ t\in [t,\bar{t});\\
\psi_j(t,B(\cdot\wedge t,\omega)-B(\bar{t},\omega)), &\mbox{if} ~t\in [\bar{t},T] ~\mbox{and}~ x(t,\omega)\in D_j.
\end{array}
\right.
\end{array}
$$
We easily know $\widetilde{u}(\cdot)\in {\cal U}[t,T]$. Therefore, we have
$$
\begin{array}{ll}
&V(t,x)\leq J(t,x;\widetilde{u}(\cdot))\\
&=\mathbb{\hat{E}}\left\{\int_t^T\Phi_1(s,x(s;t,x,\widetilde{u}(\cdot)),\widetilde{u}(s))ds+\Phi_2(x(T;t,x,\widetilde{u}(\cdot)))\right\}\\
&=\mathbb{\hat{E}}\left\{\int_t^{\bar{t}}\Phi_1(s,x(s;t,x,u(\cdot)),u(s))ds\right.\\
&\quad\left.+\mathbb{\hat{E}}\left[\int_{\bar{t}}^T\Phi_1(s,x(s;\bar{t},x(\bar{t}),u(\cdot)),\widetilde{u}(s))ds+\Phi_2(x(T;\bar{t},x(\bar{t}),\widetilde{u}(\cdot)))|{\Omega^t_{\bar{t}}}\right]\right\}\\
&=\mathbb{\hat{E}}\left\{\int_t^{\bar{t}}\Phi_1(s,x(s;t,x,u(\cdot)),u(s))ds+J(\bar{t},x(\bar{t};t,x,u(\cdot));\widetilde{u}(\cdot))\right\}\\
&\leq\mathbb{\hat{E}}\left\{\int_t^{\bar{t}}\Phi_1(s,x(s;t,x,u(\cdot)),u(s))ds+V(\bar{t},x(\bar{t};t,x,u(\cdot)))+3\varepsilon\right\},
\end{array}
$$
where the last inequality comes from (3.12). Thus, letting $\varepsilon\downarrow0$ in above inequality we get that
 $$
V(t,x)\leq \inf\limits_{u(\cdot)\in{\cal U}[t,T]}\mathbb{\hat{E}}\left\{\int_t^{\bar{t}}\Phi_1(s,x(s;t,x,u(\cdot)),u(s))ds+V(\bar{t},x(\bar{t};t,x,u(\cdot)))\right\}
$$
 by taking the infimum over $u\in {\cal U}[t,T]$. Hence the proof is complete.
\end{proof}

\vskip12pt

Let $C^{1,2}([0,T]\times{\mathbb{R}}^m; {\mathbb{R} })$ denote
 all functions ${\rm v}(t,x)$ on $[0,T]\times{\mathbb{R}}^m$ which are continuously differentiable in $t$, continuously
twice differentiable in $x$. For
${\rm v}(\cdot)\in C^{1,2}([0,T]\times{\mathbb{R}}^m;
\mathbb{R})$, define operator $L_G(u)$ by
$$
 \begin{array}{ll}
&L_G(u){\rm v}(t,x):={\rm v}_t(t, x)+ \prec{\rm v}_x(t,
x), f(t,x,u)\succ
+G(g^\top(t,x,u){\rm v}_{xx}(t, x)g(t,x,u)).
\end{array}
$$

We now give the HJB equation of the optimal
control problem 3.2 by the principle of optimality I.

\vskip12pt

{\noindent \bf Theorem 3.7} (HJB equation I) {\sl  Let Assumption 3.1 hold. If the value function of Problem 3.2 $ V(\cdot)\in C^{1,2}([0,T]\times{\mathbb{R}}^m; \mathbb{R})$, then the value function $V(\cdot)$
satisfies the HJB equation
$$
\inf_{u\in {\bf U}}\{L_G(u)V(t,x)+\Phi_1(t,x,u)\}=0,~(t,x)\in [0,T)\times{\mathbb{R}}^m\eqno(3.13)
$$
 with the terminal
condition $V(T, x)=\Phi_2(x).$
 }

\vskip12pt

\begin{proof}  From Assumption 3.1, Lemma 3.4 and Proposition 3.5, together with $u(\cdot)\in{\cal U}[t,T]$, we get
  the functions $ V(t,x)$ be bounded when $|x|$ is
bounded. Therefore, we deduce that, for each time $s\in[t,T]$,
$$
\mathbb{\hat{E}}\left[\int_t^s\prec V_x(r,x(r)),g(r,x(r),u(r))dB(r)\succ\right]=0.\eqno(3.14)
$$
 Next,  by Lemma 2.3 and (3.14) we obtain that
$$
 \begin{array}{ll}
dV(s,x(s))=&[V_t(s,x(s))+\prec V_x(s,x(s)),f(s,x(s),u(s))\succ] dt\\
&+\frac{1}{2}\prec g^\top(s,x(s),u(s))V_{xx}(s,x(s))g(s,x(s),u(s)), d<B>_s\succ\\
&+\prec V_x(s,x(s)),g(s,x(s),u(s))dB(s)\succ,
\end{array}
$$
which shows that
$$
\begin{array}{ll}
&V(s,x(s))-V(t,x)\\
&=\int_t^s[V_t(r,x(r))+\prec V_x(r,x(r)),f(r,x(r),u(r))\succ] dr\\
&\quad+\frac{1}{2}\int_t^s\prec g^\top(r,x(r),u(r))V_{xx}(r,x(r))g(r,x(r),u(r)), d<B>_r\succ\\
&\quad+\int_t^s\prec V_x(r,x(r)),g(r,x(r),u(r))dB(r)\succ.
\end{array}
\eqno(3.15)
$$
Due to Theorem 3.6, we have
  $$
V(t,x)=\inf\limits_{u(\cdot)\in{\cal U}[t,T]}\mathbb{\hat{E}}\left\{\int_t^s\Phi_1(r,x(r),u(r))dr+V(s,x(s)) \right\}.\eqno(3.16)
$$
From (3.15) and (3.16), for each $u(\cdot)\in {\cal U}[t,T]$ with $u(t)\equiv u$, by Lemma 2.4 we have
$$
\begin{array}{ll}
0&\leq \frac{1}{s-t}\mathbb{\hat{E}}\left\{\int_t^s\Phi_1(r,x(r),u(r))dr+[V(s,x(s))-V(t,x)]\right\}\\
&=\frac{1}{s-t}\mathbb{\hat{E}}\left\{\int_t^s\Phi_1(r,x(r),u(r))dr\right.\\
&\quad+\int_t^s[V_t(r,x(r))+\prec V_x(r,x(r)),f(r,x(r),u(r))\succ] dr\\
&\quad+\frac{1}{2}\int_t^s\prec g^\top(r,x(r),u(r))V_{xx}(r,x(r))g(r,x(r),u(r)), d<B>_r\succ\\
&\left.\quad+\int_t^s\prec V_x(r,x(r)),g(r,x(r),u(r))dB(r)\succ\right\}\\
&=\frac{1}{s-t}\mathbb{\hat{E}}\left\{\int_t^s\Phi_1(r,x(r),u(r))dr\right.\\
&\quad+\int_t^s[V_t(r,x(r))+\prec V_x(r,x(r)),f(r,x(r),u(r))\succ] dr\\
&\left.\quad+\frac{1}{2}\int_t^s\prec g^\top(r,x(r),u(r))V_{xx}(r,x(r))g(r,x(r),u(r)), d<B>_r\succ\right\}.
\end{array}
$$
By Assumption 3.1, (2.1) and Theorem VI-1.31 in Peng \cite{P2010}, letting $s\downarrow t$, noting $x(t)=x,u(t)=u$, we get
$$
\begin{array}{ll}
0&\leq \Phi_1(t,x, u)+V_t(t,x)+\prec V_x(t,x),f(t,x,u)\succ\\
&+\frac{1}{2}\mathbb{\hat{E}}(\prec g^\top(t,x,u)V_{xx}(t,x)g(t,x,u), v_t\succ)\\
&=\Phi_1(t,x,u)+L_G(u)V(t,x),
\end{array}
$$
which shows that
$$
0\leq \inf\limits_{u\in {\bf U}}\{\Phi_1(t,x,u)+L_G(u)V(t,x)\}.\eqno(3.17)
$$

On the other hand, for any $\varepsilon>0, 0\leq t<s\leq T$ with $s-t>0$ small enough, there exists a $u(\cdot)\equiv u_{\varepsilon,s}(\cdot)\in {\cal U}[t,T]$ such that
$$
V(t,x)+\varepsilon(s-t)\geq\mathbb{\hat{E}}\left\{\int_t^s\Phi_1(r,x(r),u(r))dr+V(s,x(s))\right\}.
$$
Thus by Lemma 2.4 we have
$$
\begin{array}{ll}
\varepsilon &\geq \frac{1}{s-t}\mathbb{\hat{E}}\left\{\int_t^s\Phi_1(r,x(r),u(r))dr+[V(s,x(s))-V(t,x)]\right\}\\
&=\frac{1}{s-t}\mathbb{\hat{E}}\left\{\int_t^s\Phi_1(r,x(r),u(r))dr\right.\\
&\quad+\int_t^s[V_t(r,x(r))+\prec V_x(r,x(r)),f(r,x(r),u(r))\succ] dr\\
&\quad+\frac{1}{2}\int_t^s\prec g^\top(r,x(r),u(r))V_{xx}(r,x(r))g(r,x(r),u(r)), d<B>_r\succ\\
&\left.\quad+\int_t^s\prec V_x(r,x(r)),g(r,x(r),u(r))dB(r)\succ\right\}\\
&=\frac{1}{s-t}\mathbb{\hat{E}}\left\{\int_t^s\Phi_1(r,x(r),u(r))dr\right.\\
&\quad+\int_t^s[V_t(r,x(r))+\prec V_x(r,x(r)),f(r,x(r),u(r))\succ] dr\\
&\left.\quad+\frac{1}{2}\int_t^s\prec g^\top(r,x(r),u(r))V_{xx}(r,x(r))g(r,x(r),u(r)), d<B>_r\succ\right\}.
\end{array}
$$

By Assumption 3.1, (2.1) and Theorem VI-1.31 in Peng \cite{P2010}, letting $s\downarrow t$, noting $x(t)=x,u(t)=u$ we get
$$
\begin{array}{ll}
\varepsilon &\geq \Phi_1(t,x, u)+V_t(t,x)+\prec V_x(t,x),f(t,x,u)\succ\\
&+\frac{1}{2}\mathbb{\hat{E}}(\prec g^\top(t,x,u)V_{xx}(t,x)g(t,x,u), v_t\succ)\\
&=\Phi_1(t,x,u)+L_G(u)V(t,x).
\end{array}
$$
Letting $\varepsilon\downarrow 0$, we have
$$
0\geq \inf\limits_{u\in {\bf U}}\{\Phi_1(t,x,u)+L_G(u)V(t,x)\}.\eqno(3.18)
$$
From (3.17) and (3.18), we get our claim. Thus the proof is complete. \end{proof}

\vskip12pt

In the rest of this section, since $-\mathbb{\hat{E}}[-\xi]=\min\limits_{P\in {\cal P}}E_P[\xi]$, we may suppose that the cost functional of our control problem is as
follows
$$
\begin{array}{rr}
\widetilde{J}(t,x; u(\cdot)):=&\min\limits_{P\in{\cal P}}E_P\left[\int_0^T\Phi_1(s, x(s), u(s))ds+\Phi_2(x(T))\right]\\
=&-\mathbb{\hat{E}}\left[-\int_0^T\Phi_1(s, x(s), u(s))ds-\Phi_2(x(T))\right].
\end{array}
\eqno(3.19)
$$

Here we notice that a controller's decision is based on a positive attitude on Knightian uncertainty.  We still consider the system of the state $x(t)$ satisfying Eq. (3.3). The optimal control problem of an objective functional (3.19) can be
stated as follows.

\vskip12pt

{\noindent \bf Problem 3.8} Select an admissible control
$\hat{u}(\cdot)\in {\cal U}[t,T]$ that minimizes $\widetilde{J}(t,x; u(\cdot))$ in (3.19) and find a value
function $\widetilde{V}$ defined by $\widetilde{V}(t,x):=\inf_{u(\cdot)\in
{\cal U}[t,T]}\widetilde{J}(t,x;u(\cdot)).$ The control $\hat{u}(\cdot)$ is called an optimal
control.

\vskip12pt

 Now define an operator $\widetilde{L}(u)$ by
$$
 \begin{array}{ll}
&\widetilde{L}(u){\rm v}(t,x):={\rm v}_t(t, x)+ \prec{\rm v}_x(t,
x), f(t,x,u)\succ
+\widetilde{G}(g^\top(t,x,u){\rm v}_{xx}(t, x)g(t,x,u)).
\end{array}
$$

\vskip12pt

{\noindent \bf Theorem 3.9} (Principle of optimality II) {\sl Let Assumption 3.1 hold. Then for any $(t,x)\in [0,T]\times {\mathbb R}^m,$ we have
$$
\begin{array}{ll}
\widetilde{V}(t,x)=\inf\limits_{u(\cdot)\in{\cal U}[t,T]}&-\mathbb{\hat{E}}\{-\int_t^{\bar{t}}\Phi_1(s,x(s;t,x,u(\cdot)),u(s))ds\\
&\quad -\widetilde{V}(\bar{t},x(\bar{t}; t,x,u(\cdot))) \},~\forall~ 0\leq t\leq\bar{t}\leq T.
\end{array}
$$
}

\vskip12pt

Next, we give the HJB equation of the optimal
control problem 3.8.

\vskip12pt

{\noindent \bf Theorem 3.10} (HJB equation II) {\sl  Let Assumption 3.1 hold. If the value function of Problem 3.8 $ \widetilde{V}(\cdot)\in C^{1,2}([0,T]\times{\mathbb{R}}^m; \mathbb{R})$, then the value function $\widetilde{V}(\cdot)$
satisfies the HJB equation
$$
\inf\limits_{u\in {\bf U}}\{\widetilde{L}(u)\widetilde{V}(t,x)+\Phi_1(t,x,u)\}=0,~(t,x)\in [0,T)\times{\mathbb~R}^m
\eqno(3.20)
$$
 with the terminal
condition $\widetilde{V}(T, x)=\Phi_2(x).$
 }

\vskip12pt

\begin{proof} First, we note that for each $u(\cdot)\in {\cal U}[t,T]$ with $u(t)\equiv u$, by Lemma 2.3 we have
$$
\begin{array}{ll}
&-\mathbb{\hat{E}}\left\{-\int_t^s\Phi_1(r,x(r),u(r))dr-[\widetilde{V}(s,x(s))-\widetilde{V}(t,x)]\right\}\\
&=-\mathbb{\hat{E}}\left\{-\int_t^s\Phi_1(r,x(r),u(r))dr\right.\\
&\quad-\int_t^s[\widetilde{V}_t(r,x(r))+\prec \widetilde{V}_x(r,x(r)),f(r,x(r),u(r))\succ] dr\\
&\quad-\frac{1}{2}\int_t^s\prec g^\top(r,x(r),u(r))\widetilde{V}_{xx}(r,x(r))g(r,x(r),u(r)), d<B>_r\succ\\
&\left.\quad-\int_t^s\prec \widetilde{V}_x(r,x(r)),g(r,x(r),u(r))dB(r)\succ\right\}\\
&=-\mathbb{\hat{E}}\left\{-\int_t^s\Phi_1(r,x(r),u(r))dr\right.\\
&\quad-\int_t^s[\widetilde{V}_t(r,x(r))+\prec \widetilde{V}_x(r,x(r)),f(r,x(r),u(r))\succ] dr\\
&\left.\quad-\frac{1}{2}\int_t^s\prec g^\top(r,x(r),u(r))\widetilde{V}_{xx}(r,x(r))g(r,x(r),u(r)), d<B>_r\succ\right\}.
\end{array}
$$
Second, the rest of proof, together with Theorem 3.9, can be obtained in a similar manner as in the proof of Theorem 3.7. Thus the proof is complete. \end{proof}

\vskip12pt

{\noindent \bf Remark 3.11} The value function for Problem 3.2 or 3.8 is not necessarily smooth enough. Hence, one uses the technique of viscosity solution to characterize the value function as the unique viscosity solution of the corresponding HJB equation (3.13) or (3.20) like that in Yong and Zhou \cite{YZ}. The related concepts of viscosity solutions are also referred to Crandall and Lions \cite{CL}, Fleming and Soner \cite{FS}, and the references therein.

\vskip12pt

\section{ Optimal consumption and portfolio policy}

 In this section, we introduce a financial market with a noise of G-Brownian motion.
 Let a sublinear expectation space  $(\Omega, {\cal H}, \{{\Omega}_t\}_{0\leq t\leq T},{\mathbb{\hat{E}}})$ host a $d$-dimensional G-Brownian motion $B(\cdot)$, where ${\Omega}_t$ denotes  an agent's available information set at instant $t$.

Let the $\{{\Omega}_t\}$-adapted processes $S_0(\cdot)$ and
$S_k(\cdot),k=1,\dots,d,$ on $[0,T]$ represent the prices of the
riskless asset and the $d$ risky assets, respectively. They satisfy
G-stochastic differential equations (G-SDEs)
 $$
 \begin{array}{ll}
dS_0(t)=&r(t)S_0(t)dt,\\
dS_k(t)=&\alpha_k(t)S_k(t)dt+\gamma^\top_k(t)S_k(t)dB(t),
\end{array}
$$
with initial prices $S_0(0)=1$ and $S_k(0)=p_k>0$. Here, $r(t)$,
$\alpha(t)=(\alpha_1(t),\cdots,\alpha_d(t))^\top$
and $\gamma_{k}(t)=(\gamma_{k1},\cdots,\gamma_{kd})^\top\ (k=1,\cdots,d)$ are deterministic
and bounded interest rate, expected returns and volatility
functions, respectively.

 We suppose that
$\gamma(t)\gamma^\top(t)$ is positive definite, where $\gamma(t)=(\gamma_1(t),\cdots,\gamma_d(t))$.
An agent chooses a portfolio
$\pi=\{\pi(t)=(\pi_1(t),\cdots,\pi_d(t))^\top, t\in[0,T]\}$,
representing the fraction of wealth invested in each risky asset. We
need a technical condition to be satisfied. A portfolio vector
process is an $\{{\Omega}_t\}$-adapted stochastic vector process
$\pi$ such that
$$\mathbb{\hat{E}}\left[\int_0^s\prec\gamma^\top(\lambda)\pi(\lambda)\pi^\top(\lambda)\gamma(\lambda),d<B>_\lambda\succ\right]<\infty \eqno(4.1)$$
for all $s\in[0,T]$. The
fraction of wealth invested in the riskless asset at time
$t\in[0,T]$ is then $1-\sum_{k=1}^d\pi_k(t)$.

Any fixed $t\in [0,T].$ Let the market price of risk be
$\theta(s):=\gamma^{-1}(s)(\alpha(s)-r(s){\bf
1})$, where ${\bf 1}=(1,\cdots,1)^\top$. Now the  evolution of the
wealth at time $s$ can be written as
$$
dX(s)=r(s)X(s)ds+X(s)\pi^\top(s)\gamma(s)[dB(s)+\theta(s)ds]-c(s)ds, ~X(t)=x,\eqno(4.2)
$$
where a consumption rate process $c=\{c(s),s\in[t,T]\}$ is a
nonnegative $\{{\Omega}_s\}$-adapted stochastic process such that
$$\mathbb{\hat{E}}\left[\int_0^sc(\lambda)d\lambda\right]<\infty,~ \forall s\in[0,T].\eqno(4.3)$$
Hence, an investor gets utility from both consumption and
wealth. Consider a problem starting at time t with
known initial condition $X(t)=x$.  Therefore, we introduce the
utility functions $U_1(\cdot)$ and $U_2(\cdot)$ of the consumption
and the wealth, respectively, which are assumed to be twice
differentiable, strictly increasing, and concave on $[0,\infty)$. The functions $I_l(\cdot),~l=1,2$ are the inverse functions of $U_l^\prime(\cdot),~l=1,2$. And,
$U^\prime_l(0)=\infty,\ U^\prime_l(\infty)=0,~ l=1,2$. It is easily
to see that there exists a constant $K>0$ such that
$$U_l(y)\leq K(1+y), ~\forall y\in [0,\infty),~l=1,2.\eqno(4.4)$$
 It is worthy of noticing that $U_l, ~l=1,2$ are not necessarily  Lipschitz continuous like function $\Phi_l, ~l=1,2$ in Section 3.  For example, the functions $U_l(y)=\frac{1}{1-\kappa} y^{\frac{1}{1-\kappa}}, ~l=1,2, 0<\kappa\neq 1$ on $[0,\infty)$ are not Lipschitz continuous.

For a pessimistic investor with multiple priors $\cal P$,  her objective function follows
$$\mathcal{J}(t,x; \pi,c)=-\mathbb{\hat{E}}\left[-\int_t^Te^{-\beta (s-t)}U_1(c(s))ds-e^{-\beta(T-t)}U_2(X(T))\right],$$
 where wealth $X(t)$ obeys (4.2), terminal
time $T>0$, $\beta$ is the utility discount rate.

Given $x\geq 0$, we say that $u=(\pi,c)$ with (4.1) and (4.3) is admissible at $(t,x)$ and write $(\pi,c)\in {\cal A}(t,x)$ if $X(s)\equiv X^{t,x,\pi,c}(s)\geq 0,~\mbox{q.s.}$ for all $s\in[t,T]$, and
$$\mathbb{\hat{E}}\left[\int_0^Te^{-\beta t}U_1^-(c(t))dt+e^{-\beta T}U^-_2(X(T))\right]<\infty.\eqno(4.5)$$
We define the
value function by
$$
\mathcal{V}(t,x)=\sup\limits_{(\pi,c)\in{\cal A}(t,x)}\mathcal{J}(t,x;\pi,c),\eqno(4.6)
$$
which shows that an agent selects consumption and investment
processes in order to maximize the sum of her expected discounted
utilities from both consumption and terminal wealth.

Now define the state price density process $p=(p_t)$ as the unique solution to following equation
$$\frac{dp_t}{p_t}=-r(t)dt-\theta^\top(t) v_t^{-1} dB(t).$$
Then $p_t$ can be used to characterize feasible consumption plans (c.f., Duffie \cite{Du}).

Set
$$
\begin{array}{ll}
{\cal X}(t,y):=&-\hat{\mathbb E}\left[-\int_t^T\zeta^t_se^{-\beta(s-t)}I_1(y\zeta_s^t)ds\right.\\
&\quad\left.-\zeta_T^te^{-\beta(s-t)}I_2(y\zeta_T^t)ds\right]<\infty,~\forall 0<y<\infty,
\end{array}
\eqno(4.7)
$$
where $\zeta_t=p_t\exp\{\int_0^t\beta(s)ds\}, \zeta_s^t=\zeta_s/\zeta_t$.
Like Lemma 4.2 of Karatzas et al. \cite{KLS}, we have that the function ${\cal X}(t,\cdot), t\in[0,T]$ defined in (4.7) is continuous and strictly decreasing on $(0,\infty)$ with ${\cal X}(t,0):=\lim_{y\downarrow 0}{\cal X}(t,y)=\infty,{\cal X}(t,\infty):=\lim_{y\rightarrow \infty}{\cal X}(t,y)=0$. Therefore for each $t\in[0,T]$ we can define by ${\cal Y}(t,\cdot):[0,\infty]\rightarrow[0,\infty]$ the inverse of the function ${\cal X}(t,\cdot)$. For any given $x>0$, define the process $C_s^{(t,x)}:=I_1({\cal Y}(t,x)\zeta_s^t),~t\leq s\leq T$ and the random variable $X_T^{(t,x)}:=I_2({\cal Y}(t,x)\zeta_T^t)$. From Karatzas et al. (\cite{KLS}, p. 1572, Eq. (6.5)), we know that, for $\forall 0< y< \infty$,
$$
\begin{array}{ll}
&\hat{\mathbb E}\left[\int_t^T\zeta^t_se^{-\beta(s-t)}I_1(y\zeta_s^t)ds+\zeta_T^te^{-\beta(T-t)}I_2(y\zeta_T^t)|{\Omega}_t\right]\\
&=-\hat{\mathbb E}\left[-\int_t^T\zeta^t_se^{-\beta(s-t)}I_1(y\zeta_s^t)ds-\zeta_T^te^{-\beta(T-t)}I_2(y\zeta_T^t)|{\Omega}_t\right].
\end{array}
\eqno(4.8)
$$
Moreover, from Lemma 2.4 and (4.8), we easily get that
$$
\begin{array}{ll}
&\hat{\mathbb E}\left[\int_t^T(\zeta^t_s)^2e^{-\beta(s-t)}I_1^\prime(y\zeta_s^t)ds+(\zeta_T^t)^2e^{-\beta(T-t)}I_2^\prime(y\zeta_T^t)|{\Omega}_t\right]\\
&=-\hat{\mathbb E}\left[-\int_t^T(\zeta^t_s)^2e^{-\beta(s-t)}I^\prime_1(y\zeta_s^t)ds-(\zeta_T^t)^2e^{-\beta(T-t)}I_2^\prime(y\zeta_T^t)|{\Omega}_t\right].
\end{array}
\eqno(4.9)
$$
Now define
$$
\begin{array}{ll}
G(t,y):=&-\hat{\mathbb E}\left[-\int_t^Te^{-\beta(s-t)}U_1(I_1(y\zeta^t_s))ds\right.\\
&\left.-e^{-\beta(T-t)}U_2(I_2(y\zeta^t_T))\right],~\forall 0<y<\infty.
\end{array}
\eqno(4.10)
$$
We give the following properties of the value function.

\vskip12pt

{\noindent\bf Proposition 4.1} {\sl Suppose that the utilities $U_l, ~l=1,2$ obey the above conditions. Then the function $G: [0,T]\times[0,\infty]\rightarrow {\mathbb R}$ defined in (4.10) is strictly decreasing, continuously differential, and satisfies
$$G_y(t,y)=y{\cal X}_y(t,y),~\forall (t,y)\in[0,T]\times(0,\infty),\eqno(4.11)$$
  $${\cal V}(t,x)=G(t,{\cal Y}(t,x)),~\forall (t,x)\in[0,T]\times(0,\infty),\eqno(4.12)$$
 $${\cal V}_x(t,x)={\cal Y}(t,x)>0,~{\cal V}_{xx}(t,x)={\cal Y}_x(t,x)<0,~ \forall (t,x)\in[0,T]\times(0,\infty).\eqno(4.13)$$
}

\vskip12pt

\begin{proof} Fix $(t,y)\in[0,T]\times (0,\infty)$. By using Lemma 2.4 and (4.8), we obtain that
$$
\begin{array}{ll}
{\cal X}(t, y)-{\cal X}(t,y-\lambda)=&-\hat{\mathbb E}\left[-\int_t^T\zeta^t_se^{-\beta(s-t)}(I_1(y\zeta^t_s)-I_1((y-\lambda)\zeta^t_s))ds\right.\\
&\left.\quad-\zeta_T^te^{-\beta(T-t)}(I_2(y\zeta_T^t)-I_2((y-\lambda)\zeta_T^t))\right]\\
&\leq-\lambda\hat{\mathbb E}\left[-\int_t^T(\zeta^t_s)^2e^{-\beta(s-t)}I_1^\prime(y\zeta^t_s)ds\right.\\
&\left.\quad-(\zeta_T^t)^2e^{-\beta(T-t)}I_2^\prime(y\zeta_T^t)\right]=:\lambda \psi(y)
\end{array}
$$
since the functions $I_1$ and $I_2$ are convex on $(0,\infty)$. Obviously, the above function $\psi(\cdot)$ is nondecreasing and continuous. Thus the left-hand derivative $D^-_y{\cal X}(t,y)$ on $y$ obeys $-\infty< D^-_y{\cal X}(t,y)\leq\psi(y).$
Similarly, we get the right-hand derivative $D^+_y{\cal X}(t,y)$ on $y$ obeys $\psi(y)\leq D^+_y{\cal X}(t,y)<\infty.$ Therefore, we shows that
$$
\begin{array}{ll}
{\cal X}_y(t, y)=&-\hat{\mathbb E}\left[-\int_t^T(\zeta^t_s)^2e^{-\beta(s-t)}I_1^\prime(y\zeta^t_s)ds\right.\\
&\left.\quad-(\zeta_T^t)^2e^{-\beta(T-t)}I_2^\prime(y\zeta_T^t)\right].
\end{array}
\eqno(4.14)
$$
On the other hand, from Lemma 4.2 of Chapter 3 in Karazatas and Shreve \cite{KS1998} and the convex of the functions $I_l,~l=1,2$ we have the inequalities
$$
\begin{array}{ll}
&\lambda(y+\lambda)(\zeta_s^t)^2I_1^\prime(y\zeta^t_s)\leq U_1(I_1(y\zeta^t_s))
-U_1(I_1((y+\lambda)\zeta^t_s))\leq\lambda y(\zeta_s^t)^2I_1^\prime((y+\lambda)\zeta^t_s),\\
&\lambda(y+\lambda)(\zeta^t_T)^2I_2^\prime(y\zeta^t_T)\leq U_2(I_2(y\zeta^t_T))-U_2(I_2(y\zeta^t_T))\leq\lambda y(\zeta^t_T)^2I_2^\prime((y+\lambda)\zeta^t_T).
\end{array}
\eqno(4.15)
$$
Thanks to (4.9), (4.10), (4.14), (4.15) and Lemma 2.4 we easily deduce that
$$\lambda(y+\lambda){\cal X}_y(t,y)\leq G(t,y+\lambda)-G(t,y)\leq \lambda y {\cal X}_y(t,y+\lambda),$$
which shows (4.11).

In a similar manner as in the discussion of (7.17) in Karazats et al. \cite{KLS}, we have that
(4.12) holds. Finally, (4.13) is easily obtained from (4.11) and (4.12). Thus the proof is complete. \end{proof}

For our aim, we now provide the following lemma.

\vskip12pt

{\noindent \bf Lemma 4.2} {\sl Let (4.1)-(4.5) hold. Then:\\
{\rm (i)} The equation (4.2) with $X(0)=x$ admits a unique solution $X(\cdot)$, and for any $T>0,\ell \geq 1$, there exists a $K_{\ell,T}>0$ such that
$$\hat{\mathbb E}\max\limits_{0\leq s\leq T}|X(s)|^\ell\leq K_{\ell,T}(1+|x|^\ell)$$
and
$$\mathbb{\hat{E}}|X(t)-X(s)|^\ell\leq K_{\ell,T}(1+|x|^\ell)|t-s|^{\ell/2}, ~\forall s,t\in[0,T].$$\\
{\rm (ii)} If $\hat{X}(0)=y\in{\mathbb R}$ and $\hat{X}(\cdot)$ is the corresponding solution of (4.2), then for any $T>0,\ell\geq 1$, there exists a $K_{\ell,T}>0$ such that

$$\mathbb{\hat{E}}\max\limits_{0\leq s\leq T}|X(s)-\hat{X}(s)|^\ell\leq K_{\ell,T}|x-y|^\ell.$$\\
{\rm (iii)} Fixed $(t,x)\in [0,T]\times(0,\infty)$. Then there exists a $\delta_x>0$ such that for each $y\in (x-\delta_x,x+\delta_x)\subset (0,\infty)$, we have
$$
\begin{array}{ll}
&|{\cal J}(t,x;u)-{\cal J}(\bar{t},y;u)|\vee|{\cal V}(t,x)-{\cal V}(\bar t, y)|\\
\leq &K_x(|x-y|+(1+x\vee y)|t-\bar{t}|^{1/2}),\quad \forall\bar{t}\in[0,T], u\in{\cal A }(t,x),
\end{array}
$$
where $K_x$ depends on $x$.
}

\vskip12pt

\begin{proof} Similar to the proof of Lemma 3.4, we can obtain (i) and (ii).

Now we let $0\leq t\leq\bar{t}\leq T$ and $x,\bar{x}\in (0,\infty)$. For any control $u\in {\cal A}(t,x)$, let $X(\cdot)$ and $\bar{X}(\cdot)$ be the states corresponding to $(t,x,u)$ and $(\bar{t},\bar{x},u)$, respectively. Therefore, by (i) and (ii) with $\ell=1$, we obtain
$$
\mathbb{\hat{E}}\sup\limits_{s\in[\bar{t},T]}|X(s)-\bar{X}(s)|\leq K_T(|x-\bar{x}|+(1+x\vee\bar{x})|t-\bar{t}|^{1/2}).
$$
Thus, for fixed $(t,x)\in[0,T]\times (0,\infty)$, due to (4.2), (4.3) and $U_l^\prime(0)=\infty,~l=1,2$ there exists $\delta_x>0$ such that
$$
\begin{array}{ll}
|{\cal J}(t,x;u)-{\cal J}(\bar{t},y;u)|\leq&  K_x (|x-y|+(1+x \vee y)|t-\bar{t}|^{1/2}),\\
&\quad\forall \bar{t}\in[0,T],~ y\in(x-\delta_x,x+\delta_x)\cap (0,\infty).
\end{array}
$$
Taking the supermum in $u\in{\cal A}(t,x)$, we obtain (iii). Thus the proof is complete.
\end{proof}

\vskip12pt

Next we provide the principle of optimality in the case of optimal consumption and portfolio which is important for deriving the equation of optimality (or HJB equation).

\vskip12pt

{\noindent \bf Theorem 4.3} (Principle of optimality III) {\sl Let (4.1)-(4.5) hold. Then for any $(t,x)\in [0,T]\times (0,\infty),$ we have
$$
\begin{array}{ll}
{\cal V}(t,x)=\sup\limits_{(\pi,c)\in{\cal A }(t,x)}&-\mathbb{\hat{E}}\{-\int_t^{\bar{t}}e^{-\beta(s-t)}U_1(c(s))ds\\
&\quad- e^{-\beta(\bar{t}-t)}{\cal V}(\bar{t},X(\bar{t})) \},~\forall ~0\leq t\leq\bar{t}\leq T.
\end{array}\eqno(4.15)
$$
}

\vskip12pt

\begin{proof} For any $\varepsilon>0$, there exists a $u=(\pi,c)\in {\cal A}(t,x)$ such that
$$
\begin{array}{ll}
&{\cal V}(t,x)-\varepsilon\leq {\cal J}(t,x;u)\\
&=-\mathbb{\hat{E}}\left\{-\int_t^Te^{-\beta(s-t)}U_1(c(s))ds-e^{-\beta(T-t)}U_2(X(T))\right\}\\
&=-\mathbb{\hat{E}}\left\{-\int_t^{\bar{t}}e^{-\beta(s-t)}U_1(c(s))ds+\mathbb{\hat{E}}\left[-\int_{\bar{t}}^Te^{-\beta(s-t)}U_1(c(s))ds-e^{-\beta(T-t)}U_2(X(T))|{\Omega^t_{\bar{t}}}\right]\right\}\\
&=-\mathbb{\hat{E}}\left\{-\int_t^{\bar{t}}e^{-\beta(s-t)}U_1(c(s))ds-e^{-\beta(\bar{t}-t)}{\cal J}(\bar{t},x(\bar{t});u)\right\}\\
&\leq-\mathbb{\hat{E}}\left\{-\int_t^{\bar{t}}e^{-\beta(s-t)}U_1(c(s))ds-e^{-\beta(\bar{t}-t)}{\cal V}(\bar{t},x(\bar{t}))\right\},
\end{array}
$$
which shows
$$
{\cal V}(t,x)\leq\sup\limits_{u\in{\cal A }(t,x)}-\mathbb{\hat{E}}\left\{-\int_t^{\bar{t}}e^{-\beta(s-t)}U_1(c(s))ds-e^{-\beta(\bar{t}-t)}{\cal V}(\bar{t},X(\bar{t})) \right\}.\eqno(4.16)
$$
On the other hand, let $a$ and $b$ satisfy $0<a<|X(s)|<b<\infty$. We
define $\tau_a£º=\inf\{s\in[0,T]:
|X(s)|\leq a\}, \tau_b£º=\inf\{s\in[0,T]:
|X(s)|\geq b\}$, and
$\tau£º=\tau_a\wedge\tau_b$, where $\tau_a,\tau_b$ and $\tau$ are $\{\Omega_s\}$-stopping time (under a sublinear expectation, the concept of stopping time is referred to Li and Peng \cite{LP}, whose operating is similar to the classical situation). $F_a^b:=[a, b]$. Obviously, $\tau\rightarrow T$ as $a\downarrow0, b\uparrow\infty$. For obtaining a necessary inequality, we cope with it by a localizing method. To this end, we now define
$$
\begin{array}{rl}
\mathcal{J}_a^b(t,x; u)&:=-\mathbb{\hat{E}}\left[-\int_t^{T\wedge\tau}e^{-\beta (s-t)}U_1(c(s))ds-e^{-\beta(T\wedge \tau-t)}U_2(X(T\wedge\tau))\right],\\
 {\cal V}_a^b(t,x)&:=\sup\limits_{u\in {\cal A}(t,x)}{\cal J}_a^b(t,x;u).
 \end{array}
 $$
 We easily know that $\mathcal{J}_a^b(t,x; u)\rightarrow\mathcal{J}(t,x; u)$, ${\cal V}_a^b(t,x)\rightarrow{\cal V}(t,x)$ as $a\downarrow0,b\uparrow\infty$.

 For any $\varepsilon>0, y\in (0,\infty)$, by Lemma 4.2, there exists a $\delta_y>0$ such that whenever $\bar{y}\in D_y=(y-\delta_y,y+\delta_y)\subset(0,\infty)$,
 $$|{\cal J}_a^b(\bar{t},y,u)-{\cal J}_a^b(\bar{t},\bar{y},u)|+|{\cal V}_a^b(\bar{t},y)-{\cal V}_a^b(\bar{t},\bar{y})|\leq\varepsilon,~\forall u\in {\cal A}(t,x), \bar{t}\in [t,T].\eqno(4.17)$$
 Since $\cup_{y\in F_a^b}D_y\supset F_a^b$, by the Borel finite covering theorem, there exists a finite open sets $D_i=(x_i-\delta_i,x_i+\delta_i), i=1,\cdots,n$ such that $\cup_{1\leq i\leq n }D_i\supset F_a^b$. For each $i$ there exists $u_i\in {\cal A}(t,x)$ such that
 $${\cal J}_a^b(\bar{t},x_i;u_i)\geq {\cal V}_a^b(\bar{t},x_i)-\varepsilon.$$
 Thus for any $x\in F_a^b$, there exists a $D_i$ such that $x\in D_i$. Hence from (4.17) we have
 $${\cal J}_a^b(\bar{t},x;u_i)\geq {\cal J}_a^b(\bar{t},x_i;u_i)-\varepsilon\geq {\cal V}_a^b(\bar{t},x_i)-2\varepsilon\geq {\cal V}_a^b(\bar{t},x)-3\varepsilon.\eqno(4.18)$$
 Then there exists a continuous function $\psi_i:[\bar{t}\wedge\tau, T\wedge\tau]\times C([0,T];{\mathbb R}^d)\rightarrow {\bf U}$ such that
 $$u_i(s,\omega)=\psi_i(s, B(\cdot\wedge s,\omega)-B(\cdot\wedge \bar{t}\wedge\tau,\omega))~\mbox{q.s.},~\forall s\in [\bar{t}\wedge\tau,T\wedge\tau].$$
 Now, given $u\in {\cal A}(t,x)$, let $X(\cdot)\equiv X(\cdot;t,x,u)$ denote the corresponding state trajectory. Define a new control
 $$
\begin{array}{ll}
\widetilde{u}(s,\omega)=\left\{\begin{array}{ll}u(s,\omega),&\mbox {if}~ s\in [t,\bar{t}\wedge\tau);\\
\psi_i(s,B(\cdot\wedge s,\omega)-B(\bar{t}\wedge\tau,\omega)), &\mbox{if} ~s\in [\bar{t}\wedge\tau,T\wedge\tau] ~\mbox{and}~ X(s,\omega)\in D_i.
\end{array}
\right.
\end{array}
$$
We easily know $\widetilde{u}=(\widetilde{\pi},\widetilde{c})\in {\cal A}(t,x)$.
Since
$$
\begin{array}{ll}
&{\cal V}_a^b(t,x)\geq \mathcal{J}_a^b(t,x; \widetilde{\pi},\widetilde{c})\\
&=-\mathbb{\hat{E}}\left\{-\int_t^{T\wedge\tau}U_1(\widetilde{c}(s))ds-U_2(X(T\wedge\tau;t,x,\widetilde{\pi},\widetilde{c}))\right\}\\
&=-\mathbb{\hat{E}}\left\{-\int_t^{\bar{t}\wedge\tau}U_1(c(s))ds+\mathbb{\hat{E}}\left[-\int_{\bar{t}\wedge\tau}^{T\wedge\tau}U_1(\widetilde{c}(s))ds-U_2(X(T\wedge\tau;t,x,\widetilde{\pi},\widetilde{c}))|{\Omega^t_{\bar{t}\wedge\tau}}\right]\right\}\\
&=-\mathbb{\hat{E}}\left\{-\int_t^{\bar{t}\wedge\tau}U_1(c(s))ds-{\cal J}_a^b(\bar{t}\wedge\tau,X(\bar{t}\wedge\tau;t,x,u);\widetilde{u})\right\}\\
&\geq-\mathbb{\hat{E}}\left\{-\int_t^{\bar{t}\wedge\tau}U_1(c(s))ds-{\cal V}_a^b(\bar{t}\wedge\tau,X(\bar{t}\wedge\tau;t,x,u))\right\}-3\varepsilon,
\end{array}
\eqno(4.19)
$$
where the last inequality comes from (4.18). Letting $a\downarrow0,b\uparrow\infty,\varepsilon\downarrow0$ in (4.19), we get that
$$
\begin{array}{ll}
&{\cal V}(t,x)\geq \sup\limits_{u\in{\cal A }(t,x)}-\mathbb{\hat{E}}\left\{-\int_t^{\bar{t}}U_1(c(s))ds-{\cal V}(\bar{t},X(\bar{t}))\right\}.
\end{array}
\eqno(4.20)
$$
Thus, combining (4.16) with (4.20) we get equation (4.15). Thus the proof is complete.
\end{proof}

\vskip12pt

Now we deduce the HJB equation
satisfied by  optimal consumption and portfolio policies. For this,
denote operator $\mathcal{L}(u)\psi$ by
$$
 \begin{array}{ll}
&\mathcal{L}(u)\psi(t,x) =\mathcal{L}(\pi,c)\psi(t,x):=\psi_t(t,x)-\beta \psi(t,x)\\
&\quad+(x\pi^\top\gamma(t)\theta(t)+r(t)x-c) \psi_{x}(t,x)+x^2\widetilde{G}(\psi_{xx}(t,x)\gamma^\top(t)\pi\pi^\top\gamma(t)).
\end{array}
$$
 Thus we give the verification theorem of
the optimal consumption and portfolio policy.

\vskip12pt

{\noindent \bf Theorem 4.4} (HJB equation III) {\sl Suppose that (4.1)-(4.5) hold. Let the
function $\psi(t,x)\in C^{1,2}([0,T]\times(0,\infty);{\mathbb{R}})$. Suppose that the function $\psi(t,x)$ satisfies the HJB equation

$$
\begin{array}{ll}
\sup\limits_{c\geq0,\pi\in {\mathbb R}^d}\left[\mathcal{L}(\pi,c)\psi(t,x)+U_1(c)\right]=0,~(t, x)\in[0,T)\times (0,\infty)
\end{array}
\eqno(4.21)
$$
 with the terminal
condition $\psi(T, x)=U_2(x).$  If
the policy $\hat{u}=(\hat{\pi},\hat{c})$ defined by
$$
\hat{u}(t)={\rm arg}\sup\limits_{c\geq0,\pi\in{\mathbb R}^d}\left[\mathcal{L}(\pi,c)\psi(t,x(t))+U_1(c(t))\right]
$$
 is admissible, then it is the  optimal consumption and
portfolio policy of the control problem (4.6). Furthermore,
$\mathcal{V}(t,x)=\psi(t,x)$.
 }

\vskip12pt
\begin{proof}  Fix $(t,x)\in[0,T)\times (0,\infty)$. Let $a$ and $b$ satisfy $0<a<|X(s)|<b<\infty$, where $X(\cdot)$ obeys (4.2). We
define the stopping times $\tau_a£º=\inf\{s\in[t,T]:
|X(s)|\leq a\}, \tau_b£º=\inf\{s\in[t,T]:
|X(s)|\geq b\}$, and
$\tau£º=\tau_a\wedge\tau_b$.

 By Lemma 4.2, we know that the function $ {\cal V}_x(t,x)$ is bounded whenever $|x|$ is
bounded. Therefore, from (4.1) we easily deduce that, for $\hat{t}\in[0,T]$,
$$
 \begin{array}{ll}
&-\mathbb{\hat{E}}\left[-\left|\int_t^{\hat{t}\wedge \tau}{\cal V}_x(s,X(s))X(s)\pi^\top(s)\gamma(s)dB(s)\right|^2\right]\\
&\quad=-\mathbb{\hat{E}}\left[-\int_t^{\hat{t}\wedge \tau}{\cal V}^2_x(s,X(s))X^2(s)\prec\gamma^\top(s)\pi(s)\pi^\top(s)\gamma(s),d<B>_s\right]<\infty.
\end{array}
$$
Thus, we get
$$
 -\mathbb{\hat{E}}\left[-\int_t^{\hat{t}\wedge \tau}e^{-\beta (s-t)}{\cal V}_x(s,X(s))X(s)\pi^\top(s)\gamma(s)dB(s)\right]=0,
$$
which shows
$$
-\mathbb{\hat{E}}\left[-\int_t^{\hat{t}}e^{-\beta (s-t)}{\cal V}_x(s,X(s))X(s)\pi^\top(s)\gamma(s)d B(s)\right]=0, \eqno(4.22)
$$
since $\tau\rightarrow T$ as $a\downarrow0,b\uparrow\infty$.

Next, we note that for each $(\pi,c)\in {\cal A}(t,x)$, by Lemma 2.3 we get that
$$
\begin{array}{ll}
&e^{-\beta\hat{t}}{\cal V}(\hat{t},X(\hat{t}))-e^{-\beta t}{\cal V}(t,x)\\
&=\int_t^{\hat{t}}e^{-\beta s}[{\cal V}_t(s,X(s))-\beta{\cal V}(s,X(s))+{\cal V}_x(s,X(s))h(s,X(s),\pi(s),c(s))] ds\\
&\quad+\frac{1}{2}\int_t^{\hat{t}}e^{-\beta s} X^2(s){\cal V}_{xx}(s,X(s))\prec\gamma^\top(s)\pi(s)\pi^\top(s)\gamma(s),d<B>_s\succ\\
&\quad+\int_t^{\hat{t}}e^{-\beta s} {\cal V}_x(s,X(s))X(s)\pi^\top(s)\gamma(s)dB(s),
\end{array}
$$
where $h(s,x,\pi,c)=(x\pi^\top\gamma(s)\theta(s)+r(s)x-c)$.
Hence from Lemma 2.4, Theorem 4.3 and (4.22) we obtain that
$$
\begin{array}{ll}
&-\mathbb{\hat{E}}\left\{-\int_t^{\hat{t}} e^{-\beta(s-t)}U_1(c(s))ds-[e^{\beta(\hat{t}-t)}{\cal V}({\hat{t}},X({\hat{t}})-{\cal V}(t,x)]\right\}\\
&=-\mathbb{\hat{E}}\left\{-\int_t^{\hat{t}} e^{-\beta(s-t)}U_1(c(s))ds\right.\\
&\quad-\int_t^{\hat{t}}e^{-\beta(s-t)}[{\cal V}_t(s,X(s))-\beta{\cal V}(s,X(s))+{\cal V}_x(s,X(s))h(s,X(s),\pi(s),c(s))] ds\\
&\quad-\frac{1}{2}\int_t^{\hat{t}} e^{-\beta(s-t)}X^2(s){\cal V}_{xx}(s,X(s))\prec\gamma^\top(s)\pi(s)\pi^\top(s)\gamma(s),d<B>_s\succ\\
&\left.\quad-\int_t^{\hat{t}}e^{-\beta(s-t)} {\cal V}_x(s,X(s))X(s)\pi^\top(s)\gamma(s)dB(s)\right\}\\
&=-\mathbb{\hat{E}}\left\{-\int_t^{\hat{t}} e^{-\beta(s-t)}U_1(c(s))dr\right.\\
&\left.\quad-\int_t^{\hat{t}}e^{-\beta(s-t)}[{\cal V}_t(s,X(s))-\beta{\cal V}(s,X(s))+{\cal V}_x(s,X(s))h(s,X(s),\pi(s),c(s))] dr\right.\\
&\quad\left.-\frac{1}{2}\int_t^{\hat{t}}e^{-\beta(s-t)} X^2(s){\cal V}_{xx}(s,X(s))\prec\gamma^\top(s)\pi(s)\pi^\top(s)\gamma(s),d<B>_s\succ\right\}.
\end{array}
$$
Hence, one the hand, we have
$$
\begin{array}{ll}
0 \geq& -\frac{1}{\hat{t}-t}\mathbb{\hat{E}}\left\{-\int_t^{\hat{t}} e^{-\beta(s-t)}U_1(c(s))dr\right.\\
&\left.\quad-\int_t^{\hat{t}}e^{-\beta(s-t)}[{\cal V}_t(s,X(s))-\beta{\cal V}(s,X(s))+{\cal V}_x(s,X(s))h(s,X(s),\pi(s),c(s))] dr\right.\\
&\quad\left.-\frac{1}{2}\int_t^{\hat{t}}e^{-\beta(s-t)} X^2(s){\cal V}_{xx}(s,X(s))\prec\gamma^\top(s)\pi(s)\pi^\top(s)\gamma(s),d<B>_s\succ\right\}.
\end{array}
$$
Consequently, by Theorem VI-1.31 in Peng \cite{P2010} and (2.1), letting $\hat{t}\downarrow t$ we deduce
$$
0\geq \sup\limits_{c\geq0,\pi\in {\mathbb R}^d}\{U_1(c)+{\cal L}(\pi,c){\cal V}(t,x)\}.\eqno(4.23)
$$
On the other hand, for any $\varepsilon>0, 0\leq t<\hat{t}<T$ with $\hat{t}-t>0$ small enough, there exists a $(\pi,c)\in{\cal A}(t,x)$ such that
$$
\begin{array}{ll}
-\varepsilon \leq & -\frac{1}{\hat{t}-t}\mathbb{\hat{E}}\left\{-\int_t^{\hat{t}} e^{-\beta(s-t)}U_1(c(s))dr\right.\\
&\left.\quad-\int_t^{\hat{t}}e^{-\beta(s-t)}[{\cal V}_t(s,X(s))-\beta{\cal V}(s,X(s))+{\cal V}_x(s,X(s))h(s,X(s),\pi(s),c(s))] dr\right.\\
&\quad\left.-\frac{1}{2}\int_t^{\hat{t}}e^{-\beta(s-t)} X^2(s){\cal V}_{xx}(s,X(s))\prec\gamma^\top(s)\pi(s)\pi^\top(s)\gamma(s),d<B>_s\succ\right\}.
\end{array}
$$
Thus, by Theorem VI-1.31 in Peng \cite{P2010} and (2.1), Letting $\hat{t}\downarrow t, \varepsilon\downarrow 0$, we get that
$$
0\leq \sup\limits_{c\geq0,\pi\in {\mathbb R}^d}(U_1(c)+{\cal L}(\pi,c){\cal V}(t,x)\}.\eqno(4.24)
$$
Thanks to (4.23) and (4.24), HJB equation (4.21) holds. By Proposition 4.1, we know that ${\cal V}\in C^{1,2}([0,T]\times(0,\infty); {\mathbb R})$, which shows our assertion.  Thus the proof is complete. \end{proof}

\vskip12pt

In what follows, we utilize Theorem 4.4 to characterize the optimal consumption and portfolio policy. Since the consumption policy can be solved by the concave
maximization in (4.21), we get the optimal feedback consumption
policy
$$
\begin{array}{ll}
\hat{c}(t,x) =I_1(\mathcal{V}_x(t,x)).
\end{array}
\eqno(4.25)
$$

 For the portfolio policy, we make a quadratic maximization problem in (4.21). Due to Proposition 4.1, we know that ${\cal V}_{xx}(t,x)<0$. Therefore, there exists a positive definite matrix $\bar{\Lambda}\in \Sigma$ such that
 $$
 \begin{array}{ll}
  \widetilde{G}({\cal V}_{xx}(t,x)\gamma^\top(t)\pi(t)\pi^\top(t)\gamma(t))&=\frac{1}{2}{\cal V}_{xx}(t,x)tr({\bar{\Lambda}}\gamma^\top(t)\pi(t)\pi^\top(t)\gamma(t))\\
 &=\frac{1}{2}{\cal V}_{xx}(t,x)\pi^\top(t)\gamma(t)\bar{\Lambda}\gamma^\top(t)\pi(t),
 \end{array}
  $$
 which shows that
 $$
 \begin{array}{ll}
 \frac{\partial}{\partial \pi}\widetilde{G}({\cal V}_{xx}(t,x)\gamma^\top(t)\pi(t)\pi^\top(t)\gamma(t))={\cal V}_{xx}(t,x)\gamma(t)\bar{\Lambda}\gamma^\top(t)\pi(t).
  \end{array}
 $$
Hence, from HJB equation (4.21) we get the optimal portfolio as follows
$$
\hat{\pi}(t,x)=-\frac{(\gamma^\top(t))^{-1}{\bar{\Lambda}^{-1}}\theta(t)\mathcal{V}_x(t,x)}{x\mathcal{V}_{xx}(t,x)}.
\eqno(4.26)
$$

Now we state the following modification of the classical Mutual Fund
Theorem.

\vskip12pt

{\noindent\bf Theorem 4.5} {\sl In the above financial market, we
have: \\
{\rm (i)} The optimal portfolio involves an allocation between the
risk-free fund $F^1(t)$ and a risky fund that consists only of risky
assets:
$F^2(t)=(\gamma^\top(t))^{-1}{\bar{\Lambda}}^{-1}\theta(t)$, where
the vector $F^2(t)$ represents the portfolio
weights of the risky
assets at time $t$. \\
{\rm (ii)} The optimal proportional allocations
$\varpi^k(t,x)$ of wealth in the fund $F^k(t), k
= 1, 2$, at time $t$ are given by
$$
\begin{array}{ll}
\varpi^2(t,x)&=-\frac{\mathcal{V}_{x}(t,x)}{x\mathcal{V}_{xx}(t,x)},\\
\varpi^1(t,x)&=1-\varpi^2(t,x).
\end{array}
 $$
 }

\vskip12pt
\begin{proof} We know that the right-hand side of (4.26) equals $\varpi^2(t,x) F^2(t)$.
  Hence, the proof is complete. \qquad\end{proof}

\vskip12pt
In the formula (4.26), since $\Sigma$ reflects an investor's uncertainty on the
financial environment which causes her multiple-priors ${\cal P}$, the investor's optimal consumption and portfolio decision with pessimism is based on the factor $\bar{\Lambda}$. In order to show the effect of an investor's Knightian uncertainty on optimal portfolio policy, we will give an illustrative example.

\vskip12pt

{\noindent\bf Example 4.6} Optimal policies for CRRA utility.

\vskip12pt

Here, we discuss a special case.  Let the utility
functions of an investor be
$$U_l(z)=\frac{z^{1-\kappa}}{1-\kappa},\ \kappa>0,\ \kappa\neq 1,\ l=1,2.$$
Thus, we obtain $I_l(y)=y^{-1/\kappa}$. From (4.21) of
Theorem 4.4 and optimal policy $(\pi,c)$ of (4.25) and (4,26),  we
easily deduce that the value function $\mathcal{V}$ satisfies the HJB equation
as follows
$$
\begin{array}{ll}
&\mathcal{V}_t(t,x)-\beta
\mathcal{V}(t,x)+r(t)x\mathcal{V}_x(t,x)+\frac{\kappa}{1-\kappa}(\mathcal{V}_x(t,x))^{\frac{\kappa-1}{\kappa}}\\
&\quad-\frac{{\cal V}_x^2(t,x)}{2\mathcal{V}_{xx}(t,x)}\theta^\top(t){\bar{\Lambda}}^{-1}\theta(t)=0
\end{array}
\eqno(4.27)
$$
with the terminal condition
$\mathcal{V}(T,x)=U_2(x)=\frac{1}{1-\kappa}x^{1-\kappa}$.

We guess the form of solution to the equation (4.27) as follows
$$\mathcal{V}(t,x)=\frac{1}{1-\kappa}A^\kappa(t)x^{1-\kappa},$$
which shows that
$$
 \kappa A^{\kappa-1}(t)A^\prime(t)-\eta(t) A^\kappa(t)+
\kappa A^{\kappa-1}(t)=0, ~ A(T)=1,\eqno(4.28)
$$
where
$$
\eta(t):=\beta-(1-\kappa)r(t)-\frac{1-\kappa}{2\kappa}\theta^\top(t){\bar{\Lambda}}\theta(t).
$$

We easily know that equation (4.28) has
the unique solution $A(t).$  In terms of (4.25) and (4.26), we get that
$$
\begin{array}{ll}
\hat{c}(t,x)=&\frac{x}{A(t)},\\
\hat{\pi}(t,x)=&\frac{1}{\kappa}\left(\gamma^\top(t)\right)^{-1}{\bar{\Lambda}}^{-1}\theta(t).
\end{array}
$$

Especially, if we only consider a risky asset (i.e. $d=1$), then we obtain $\bar{\Lambda}=\bar{\sigma}^2$. Suppose that the parameters $r(t)=r,\alpha(t)=\alpha,\gamma(t)=\gamma$ are constants, which shows that $\eta(t)=\eta$ is a constant. Then we have the solution to equation (4.28) as follows
$$
A(t)=\left(\frac{\kappa}{\eta}+(1-\frac{\kappa}{\eta})e^{-\frac{\eta}{\kappa}(T-t)}\right)^{-1}.\eqno(4.29)$$
Thus we have the following corollary.

\vskip12pt

{\noindent \bf Corollary 4.7} {\sl Let the riskless rate, the mean return rate and the volatility of the risky asset be the constants $r,\alpha$ and $\gamma$, respectively. Then the optimal consumption and portfolio policy is as follows
$$
\begin{array}{ll}
\hat{c}(t,x)=&\frac{x}{A(t)},\\
\hat{\pi}(t,x)=&\frac{1}{\kappa\gamma^2(t)\bar{\sigma}^2}(\alpha(t)-r(t)),
\end{array}
\eqno(4.30)
$$
where $A(t)$ is given by (4.29).
}

 The formula (4.30) shows that an investor's consumption rate is
higher as she is wealthier. Also, the higher the market price of risk and the lower the
risk of the risky asset is, the more invests in the
risky asset. An investor with the higher Knightian uncertainty degree puts less in the risky asset. Our result which depends on the investor's Knightian uncertainty of the market volatility is somewhat consistent
with the results in \cite{CE,FeiINS,FeiSM}.

Finally, in the rest of this section, we explore the optimal consumption and portfolio for an optimistic investor,  her objective function follows
$${\mathbb J}(t,x; \pi,c)=\mathbb{\hat{E}}\left[\int_t^Te^{-\beta (s-t)}U_1(c(s))ds+e^{-\beta(T-t)}U_2(X(T))\right],$$
 where wealth $X(s)$ with $X(t)=x$ follows Eq. (4.2).

Given $x\geq 0$, we say that $u=(\pi,c)$ is admissible at $(t,x)$ and write $(\pi,c)\in {\cal A}(t,x)$ if $X(s)\equiv X^{t,x,\pi,c}(s)\geq 0~\mbox{q.s.},$ for all $s\in[t,T]$.

We now define the value function by
$$
{\mathbb V}(t,x)=\sup\limits_{(\pi,c)\in{\cal A}(t,x)}{\mathbb J}(t,x;\pi,c).\eqno(4.31)
$$

For our aim, we provide the following principle of optimality which is an analogue of Theorem 4.3, its proof is omitted.

\vskip12pt

{\noindent \bf Theorem 4.8} (Principle of optimality VI) {\sl Let (4.1)-(4.5) hold. Then for any $(t,x)\in [0,T]\times (0,\infty),$ we have
$$
\begin{array}{ll}
{\mathbb V}(t,x)=\sup\limits_{(\pi,c)\in{\cal A }(t,x)}&\mathbb{\hat{E}}\{\int_t^{\bar{t}}e^{-\beta(s-t)}U_1(c(s))ds\\
&\quad+e^{-\beta(\bar{t}-t)}{\mathbb V}(\bar{t},X(\bar{t})) \},~\forall ~0\leq t\leq\bar{t}\leq T.
\end{array}
$$
}

Next, we derive the HJB equation
satisfied by  optimal consumption and portfolio policies. For this,
denote operator ${\mathbb L}_G(u)\psi$ by
$$
 \begin{array}{ll}
&{\mathbb L}_G(u)\psi(t,x) ={\mathbb L}_G(\pi,c)\psi(t,x)\\
&:=\psi_t(t,x)-\beta \psi(t,x)+
(x\pi^\top\gamma(t)\theta(t)+r(t)x-c) \psi_{x}+x^2G(\psi_{xx}(t,x)\gamma^\top(t)\pi\pi^\top\gamma(t)).
\end{array}
$$
 We now give the criterion of
the optimal policy for an optimistic investor. The following result can be similarly proven as in the proof of Theorem 4.4.

\vskip12pt

{\noindent \bf Theorem 4.9} (HJB equation VI) {\sl Let (4.1)-(4.5) hold. Let the
function $\psi(t,x)\in C^{1,2}([0,T]\times(0,\infty); {\mathbb{R}})$. Suppose that the function $\psi(t,x)$ satisfies the HJB equation

$$
\begin{array}{ll}
\sup\limits_{c\geq0,\pi\in {\mathbb R}^d}\left[{\mathbb L}_G(\pi,c)\psi(t,x)+U_1(c)\right]=0,~(t,x)\in[0,T)\times(0,\infty)
\end{array}
$$
 with the terminal
condition $\psi(T, x)=U_2(x).$  If
the policy $\hat{u}=(\hat{\pi},\hat{c})$ defined by
$$
\hat{u}(t)={\rm arg}\sup\limits_{c\geq0,\pi\in {\mathbb R}^d}\left[{\mathbb L}_G(\pi,c)\psi(t,x(t))+U_1(c(t))\right]
$$
 is admissible, then it is the  optimal consumption and
portfolio policy of the control problem (4.31). Furthermore,
${\mathbb V}(t,x)=\psi(t,x)$.
 }

\vskip12pt

\section{Conclusion}

Firstly, after the calculus of a G-Brownian motion and a sublinear expectation is introduced, in Section 3, we set up the optimum principles of stochastic controls with ambiguity while HJB equations for our optimal controls are derived. Here we suppose the value function of the optimal stochastic control problem is sufficiently smooth such that HJB equation admits a classical solution. However, if the value function is insufficiently smooth, then the corresponding HJB equation only has a viscosity solution. To this end, we will further discuss our optimal stochastic control and the corresponding viscosity solution of HJB equation under weaker conditions.

Secondly, in Section 4 the optimal consumption and portfolio is characterized through providing the optimum principle and the related HJB equation for maximizing agent's expectation utility with ambiguity. For the utility imposed some conditions, an agent's value function is enough smooth such that the corresponding HJB equation admits a classical solution. Hence a verification theorem of an optimal decision is proven, and then the modified two-fund separation theorem is given. For an aim of explanation, we provide an illustrative example which shows that for a pessimistic investor with ambiguity, the less an investor invests in a risky asset, the more uncertain she is.

  Next, we further compare our results with the ones from a classical optimal stochastic control. In Yong and Zhou \cite{YZ}, they consider an optimal stochastic control with an unambiguous environment, and provide the optimum principles and the corresponding HJB equations. However, our results are different from theirs since the term with uncertainty in our results appears. Especially, our financial application of the optimal consumption and portfolio can provide a new perspective for the study of the behaviour finance as an agent is uncertain in a complex environment.

  Finally, since it is doubtful that a real system is disturbed by a classical Brownian motion, the notion of a fractional Brownian motion with the self-similarity and the long range dependence is put forward (c.f., Mandelbrot and Van Ness \cite{MN}).
   Recently, Chen \cite {Chen} studied a fractional G-white noise theory, wavelet decomposition for fractional G-Brownian motion, and bid-ask pricing application to finance under Knightian uncertainty. Hence our question is how the stochastic differential equations driven by a fractional G-Brownian motion with volatility ambiguity is characterized. Moreover, we will study the optimal stochastic control problem with an ambiguous fractional G-Brownian motion in the future.

\vskip12pt

\section*{Acknowledgments} The work is supported by National Natural Science Foundation of China (71171003, 71271003, 71210107026), Anhui Natural Science Foundation
 (10040606003), and Anhui Natural
Science Foundation of Universities (KJ2012B019, KJ2013B023).

\end{document}